\documentclass{article}
\usepackage{amsthm,amsmath,amsfonts,amssymb}

\newlength{\hchng}
\newlength{\vchng}
\newtheorem{thm}{Theorem}[section]
\newtheorem{prop}[thm]{Proposition}

\newtheorem{lemma}[thm]{Lemma}
\newtheorem{defn}[thm]{Definition}

\newtheorem{preremark}[thm]{Remark}
\newenvironment{remark}{\begin{preremark}\rm}{\medskip \end{preremark}}
\numberwithin{equation}{section}

\newcommand{\R}{\mathbb R}

\DeclareMathOperator{\dv}{div}

\begin{document}

\title{On  possible isolated blow-up phenomena of the 3D-Navier-Stokes equation
and a regularity criterion in terms of supercritical function space condition and smoothness condition along the streamlines}
\author{Chi Hin Chan and Tsuyoshi Yoneda}
\maketitle
\begin{center}
Institute for Mathematics and its Applications, University of Minnesota\\
207 Church Street SE, Minneapolis, MN 55455-0134, USA\\
\end{center}
\begin{center}
Department of Mathematics and Statistics,
University of Victoria \\
PO Box 3060 STN CSC,
Victoria, BC, Canada, V8W 3R4 \\
\end{center}

\bibliographystyle{plain}
\noindent{\bf Abstract:} The first goal of our paper is to give a new type of regularity criterion
for solutions $u$ to Navier-Stokes equation in terms of some supercritical function space condition $u \in L^{\infty}(L^{\alpha ,*})$ 
(with $\frac{3}{4}(17^{\frac{1}{2}} -1)< \alpha < 3  $)
 and some exponential control on the growth rate of $\dv (\frac{u}{|u|})$ along the streamlines of u. This regularity
criterion greatly improves the previous one in \cite{smoothnesscriteria}. The proof
leading to the regularity criterion of our paper basically follows the one
in \cite{smoothnesscriteria}. However, we also point out that totally new idea which
involves the use of the new supercritical function space condition is necessary
for the success of our new regularity criterion in this paper.

The second goal of our paper is to construct a divergence free vector field $u$
within a flow-invariant tubular region with increasing twisting of streamlines
towards one end of a bundle of streamlines. The increasing twisting of streamlines
is controlled in such a way that the associated quantities
$\|u\|_{L^{\alpha}}$ and $\|\dv(\frac{u}{|u|})\|_{L^{6}}$ blow up while preserving the
finite energy property $u\in L^{2}$ at the same time. The purpose of such a
construction is to demonstrate the necessity to go beyond the scope
covered by some previous regularity criteria  such as \cite{Vasseur2} or
\cite{Esca}.
We also briefly mention how this construction is related to the regularity
criterion proved in our paper.



\section{Introduction}

The first goal of this paper is to give a new type of regularity criterion of solutions $u$ to the Navier-Stokes equation in terms of 
some weak $L^{\alpha}$ space condition on the velocity $u$ (with some $\frac{3}{4} (17^{\frac{1}{2}} -1)<\alpha<3$  ) and some
exponential control of $\dv{\frac{u}{|u|}}$ along the streamlines. 
The second goal of this paper is to give \emph{possible} blow-up situations for 3D-Navier-Stokes equation through the construction of a finite energy divergence free
velocity field $u$ with $u \notin L^{\alpha}$ with $2< \alpha <3$ and $\dv u \notin L^{6}$.
The Navier-Stokes equation on
$\mathbb{R}^3$ is given by

\begin{equation}\label{NS}
\begin{cases}
\partial_{t} u -\triangle u + \dv (u\otimes u) + \nabla P  = 0,  \\
div (u)  = 0,\quad u|_{t=0}=u_0
\end{cases}
\end{equation}
in which 
$u$ is a vector-valued function representing the velocity of
the fluid, and $P$ is the pressure. The  initial  value problem of the above
equation is endowed with the condition that $u(0, \cdot ) = u_{0}
\in L^2(\mathbb{R}^3)$. 

Modern regularity theory for solutions to equation \eqref{NS} began with the works of Leray \cite{Leray} and Hopf \cite{Hopf} in which they established, with respect to
any given initial datum $u_{0} \in L^{2}(\mathbb{R}^{3})$ which is weakly divergence free, 
the existence of a weak solutions $u : [0 , \infty ) \times \mathbb{R}^{3} \rightarrow \mathbb{R}^{3} $ lying in the class of 
$L^{\infty}(0, \infty ;
L^2(\mathbb{R}^3)) \cap L^2(0, \infty ; \dot{H}^1(\mathbb{R}^3))$ which satisfies the global energy inequality. Since the time of Leary and Hopf, any weak solution 
to equation \eqref{NS} which satisfies the finite energy, finite dissipation, and global energy inequalities is called Leray-Hopf solutions to \eqref{NS}.\

After the fundamential works of Leray and Hopf, progress in addressing the full regularity of Leray-Hopf solutions has been very slow. It was only in 1960 that 
significant progress was made by Prodi \cite{Prodi}, Serrin \cite{Serrin}, Ladyzhenskaya \cite{Ladyzhenskaya}, and their joint efforts lead to 
the following famous Prodi-Serrin-Ladyzhenskaya criterion for Leray-Hopf solutions (see the introduction of \cite{Esca} for more historical remarks about this).   
 
\begin{thm}\label{heatequationpertubation}
[Prodi, Serrin, Ladyzhenskaya] Let $u \in L^{\infty}(0, T ;
L^{2}(\mathbb{R}^{3}) ) \cap L^{2}(0, T ;\dot{H}^{1}(\mathbb{R}^{3}))$ be a
Leray-Hopf weak solution to \eqref{NS}, which also
satisfies $u \in L^{p}(0,\infty ;
L^{q}(\mathbb{R}^{3}))$, for some $p, q$ satisfying $\frac{2}{p} + \frac{3}{q} = 1,$
with $q > 3$. Then, $u$ is smooth on $(0,T]\times \mathbb{R}^{3}$ and is
uniquely determined in the following sense
\begin{itemize}
\item suppose $v \in L^{\infty}(0, T ; L^{2}(\mathbb{R}^{3}) ) \cap L^{2}(0, T ;
\dot{H}^{1}(\mathbb{R}^{3}))$ is another Leray-Hopf weak solution such that $u(0, \cdot ) = v(0,\cdot )$. Then, it follows that
$u=v$ on $(0,T]\times \mathbb{R}^{3}$. 
\end{itemize}
\end{thm}

The success of the Prodi-Serrin-Ladyzhenskaya criterion was based on the fact that the integral condition $u \in L^{p}(L^{q})$
with $p$, $q$ satisfying $\frac{2}{p} + \frac{3}{q} = 1$ \emph{and} $q > 3$ ensures that the Leray-Hopf solution $u$ behaves \emph{like} a solution to a slightly 
pertubated heat equation. It is also worthwhile to mention that the exceptional case of $u \in L^{\infty}(L^{3})$ was missed in the above regularity citerion of Prodi,
Serrin, and Ladyzhenskaya, and it was not until very recently that the regularity of solutions in the exceptional case $u \in L^{\infty}(L^{3})$ was finally established 
in the famous work \cite{Esca} due to L. Escauriaza, G. Seregin, and V. Sverak.\

After the appearance of the Prodi-Serrin-Ladyzhenskaya criterion, many different regularity cirteria of solutions to \eqref{NS} was established by researchers working in
the regularity theory of \eqref{NS}. Among these, for instance, Beir\~ao da Veiga established in \cite{Veiga} a regularity criterion in terms of the integral condition 
$\nabla u \in L^{p}(0,\infty ; (L^{q}(\mathbb{R}^{3})))$ with $\frac{2}{p} + \frac{3}{q} = 2$ (and $1 < p < \infty$)  imposed on $\nabla u$. In the same spirit of  
\cite{Veiga},  Beale, Kato and Majda \cite{Beale} gave a regularity criterion for solutions $u$ to \eqref{NS} in terms of the condition 
$\omega \in L^1(0,\infty;L^\infty(\R^3))$ imposed on the vorticity $\omega=\mathrm{curl}\ u$ associated to $u$. This regularity criterion was further improved by Kozono 
and Taniuchi in \cite{Kozono} (see also \cite{Nakai-Yoneda}). Besides these, other important works such as \cite{ChenStrainTsaiYauII} and \cite{TypeIblowupKNSS}, in which
type I blow up was excluded for axisymmetric solutions to \eqref{NS}, are attracting a lot of attentions. Due to the limitation of space and the vast 
literature in the regularity theory for solutions to \eqref{NS}, we do not try to do a complete survey here.\

However, we would like to mention an interesting regularity criterion in \cite{Vasseur2} due to Vasseur, since it is related to the main result of this paper and also to 
the previous partial result \cite{smoothnesscriteria} by the first author. \cite{Vasseur2} gave a regularity criterion for solutions $u$ to \eqref{NS} in terms of the 
integral condtion $\dv(\frac{u}{|u|}) \in L^{p}(0,\infty ;L^{q}(\mathbb{R}^{3}))$ with $\frac{2}{p}+ \frac{3}{q} \leqslant \frac{1}{2}$ 
imposed on the scalar quantity $F = \dv (\frac{u}{|u|})$.

One of the main purposes of this paper, however, is to establish the following regularity criterion for solutions $u$ to \eqref{NS} in terms of some exponential 
control on the rate of change of $F = \dv (\frac{u}{|u|})$ along the streamlines of $u$ and some \emph{weak $L^{\alpha}$ space} condition imposed on $u$. 

\begin{thm}\label{goal}

 Let $u \in L^{\infty}(0,T ; L^{2}(\mathbb{R}^{3}))\cap L^{2}(0,T; \dot{H}^{1}(\mathbb{R}^{3}))$ be a Leray Hopf solution 
to \eqref{NS} which is smooth up to a possible blow up time $T$ with $u_0\in \mathcal
S(\mathbb{R}^3)$. Let us assume that $u$ and $F = \dv (\frac{u}{|u|})$ satisfy the following  conditions.

\begin{itemize}
\item $u \in L^{\infty}(0,T;L^{\alpha ,\infty}(\mathbb{R}^{3}))$, for some given $\alpha \in (2,3)$ which satisfies $1+ 2(\frac{\alpha}{3} + \frac{3}{\alpha}) > 0$.
\item There exists some $r_{0} > 0$ and $M_{0} > 0$ such that $|u| \leqslant M_{0}$ is valid on the region 
$[0,T)\times \{x\in \mathbb{R}^{3} : |x| \geqslant r_{0}\}$

\item
For some given constants $A > 0$ and $L > 0$, the property $|\frac{u\cdot \nabla F}{|u|}|\leq A|F|$ is valid on 
$\{(t,x) \in [0,T)\times B(r_{0}) : |F(t,x)| \geqslant L\} $ (Here, $B(r_{0}) = \{x \in \mathbb{R}^{3} : |x| < r_{0}\}$). 
\end{itemize}

Then the smoothness of $u$ can be
extended beyond the time $T$.
\end{thm}

Here, we give a few remarks which illustrate the significance of Theorem \ref{goal}. We start with the third condition in Theorem \ref{goal} in which we see the condition
$|\frac{u\cdot \nabla F}{|u|}|\leq A|F|$ imposed on the region $[0,T)\times \{x\in \mathbb{R}^{3} : |x| \geqslant r_{0}\}\cap \{|F| \geqslant L\}$.  
We can see the geometric 
meaning of the constraint $|\frac{u\cdot \nabla F}{|u|}|\leq A|F|$ on $[0,T)\times \{x\in \mathbb{R}^{3} : |x| \geqslant r_{0}\}\cap \{|F| \geqslant L\}$ if we recast 
it in the following geometric language.

\begin{itemize}
\item For any time slice $t \in [0,T)$, and any streamline $\gamma : [0, S) \rightarrow \mathbb{R}^{3}$ of the veclocity profile 
$u(t,\cdot )$ which is \emph{parameterized by arclength} (that is, $\frac{d\gamma}{ds} = \frac{u}{|u|}(\gamma (s))$) and with image $\gamma ([0,S))$ lying in the 
region $ \{x\in \mathbb{R}^{3} : |x| \geqslant r_{0} , |F(t,x)| \geqslant L\}$, we have $|\frac{d}{ds} (F(\gamma (s)))| \leqslant A\cdot |F(\gamma (s))|$, for any 
$0\leqslant s \leqslant S$.
\end{itemize}

The condition $|\frac{d}{ds} (F(\gamma (s)))| \leqslant A\cdot |F(\gamma (s))|$ gives some exponential control on $F$ \emph{along each streamline} of the fluid
within the space region on which both $u$ and $F = \dv (\frac{u}{|u|})$ are large. Our \emph{original motivation} \emph{was} to prove that the smoothness of the solution 
$u : [0,T)\times \mathbb{R}^{3} \rightarrow \mathbb{R}^{3}$ to \eqref{NS} can be extended beyond the possible blow up time $T$ under the third condition of Theorem 
\ref{goal} and the Leray-Hopf property $u \in L^{\infty}(L^{2})\cap L^{2}(\dot{H}^{1})$ of the solution. But our experience told us that this cannot be so easily 
achieved \emph{without} the involvement of the following additional condition (which is the first condition of Theorem \ref{goal} ). 

\begin{itemize}
\item $u \in L^{\infty}(0,T:L^{\alpha ,\infty}(\mathbb{R}^{3}))$, for some given $\alpha \in (2,3)$ with $1+ 2(\frac{\alpha}{3} + \frac{3}{\alpha}) > 0$. 
\end{itemize}

To clarify the \emph{necessity} of the condition $u \in L^{\infty}(L^{\alpha ,\infty})$ with some $\alpha \in (2,3)$ satsifying 
$1+ 2(\frac{\alpha}{3} + \frac{3}{\alpha}) > 0$, let us mention a piece of work \cite{smoothnesscriteria} by the first author in which smoothness of a Leray-Hopf solution
$u : [0,T) \times \mathbb{R}^{3} \rightarrow \mathbb{R}^{3}$ is established beyond the possible blow up time $T$ under the following condtion. 

\begin{itemize}
\item (condition in the regularity criterion of \cite{smoothnesscriteria}) $|\frac{u\cdot \nabla F}{|u|^{\delta}}| \leqslant A |F|$ is valid on 
$[0,T)\times \mathbb{R}^{3}$, with $A > 0$ to be a given constant and $\delta$ to be a given constant with $0 < \delta < \frac{1}{3}$. 
\end{itemize}

The above mentioned regularity criterion based on the condition  $|\frac{u\cdot \nabla F}{|u|^{\delta}}| \leqslant A |F|$ (with $0 < \delta < \frac{1}{3}$) was
established in \cite{smoothnesscriteria} through applying the De Giorgi's method as developed by A. Vasseur in \cite{Vasseur}.
The main idea of the De Giorgi's method in \cite{Vasseur} is based on the establishment of the following nonlinear recurrence relation of the energy $U_{k}$ of 
a truncated function $v_{k} = [|u| - R(1-\frac{1}{2^{k}})]_{+}$ of the solution $u$ to \eqref{NS} over a certain space time region (for a precise definition of $U_{k}$,
see section 3 of this paper, or alternatively \cite{Vasseur} or \cite{smoothnesscriteria}).
\begin{equation}\label{nonlinearrecurrence}
U_{k} \leq \frac{C_{0}^{k}}{R^{\lambda}} U_{k-1}^{\beta} 
\end{equation}
According to the idea in \cite{Vasseur}, for a given solution $u$ to \eqref{NS} on $[0,T)\times \mathbb{R}^{3}$ with possible blow up time $T$, 
the $L^{\infty}$-boundedness \emph{conclusion} $|u| \leqslant R$ over $[\frac{T}{2},T)\times \mathbb{R}^{3}$ (for some
sufficiently large $R$) can be drawn from relation \eqref{nonlinearrecurrence} provided one can ensure that $\beta > 1$ and $\lambda > 0$ are valid 
\emph{simultaneously}. Roughly speaking, $\lambda > 0$ ensures the smallness of the energy $U_{1}$ of the first truncated function $v_{1}$, due to the fact that 
$\frac{1}{R^{\lambda}}$ will become small as $R$ is sufficiently large. The smallness of $U_{1}$ will tragger the nonlinear recurrence effect of relation 
\eqref{nonlinearrecurrence} which eventually causes the very fast decay of $U_{k}$ to $0$ (see Lemma \ref{Vass} which originally appeared in \cite{Vasseur} ).
This resulting decay of $U_{k}$ to $0$ then implies the desired boundedness conclusion $|u| \leqslant R$ over $[\frac{T}{2},T)\times \mathbb{R}^{3}$, which in turn
extends the smoothness of $u$ beyond the possible blow up time $T$. However, it was illustrated in \cite{smoothnesscriteria} that the requirement that $\beta > 1$ \emph{and}    
$\lambda > 0$ have to hold simultaneously \emph{prevents} us to push the constant $\delta$ (in the condition $|\frac{u\cdot \nabla F}{|u|^{\delta}}| \leqslant A |F|$) 
to go beyond the range $(0 , \frac{1}{3})$. This limitation of the De Giorgi method of \cite{Vasseur} basically comes 
from the fact that the index $\beta$ in relation \eqref{nonlinearrecurrence} is typically $\frac{5}{3}$ or $\frac{4}{3}$, which is \emph{too large} for the survival of
the condition $\lambda > 0$ in the same relation \eqref{nonlinearrecurrence}.\

As a result, the use of the extra condition  
$u \in L^{\infty}(L^{\alpha ,\infty})$, with $\alpha \in (2,3)$ satisfying $1+ 2(\frac{\alpha}{3} + \frac{3}{\alpha}) > 0$ can help us to 
\emph{lower} the index $\beta$ of relation \eqref{nonlinearrecurrence} from the typical $\frac{5}{3}$ or $\frac{4}{3}$ to become \emph{as close to $1$ as possible}, and 
this in turn ensures the survival of $\lambda > 0$ in the same relation \eqref{nonlinearrecurrence}. On the other hand, we have to address the question of whether 
the condition of $u \in L^{\infty}(L^{\alpha ,\infty})$, with $\alpha \in (2,3)$ satisfying $1+ 2(\frac{\alpha}{3} + \frac{3}{\alpha}) > 0$ is too strong as an 
assumption. Note that the constraint $1+ 2(\frac{\alpha}{3} + \frac{3}{\alpha}) > 0$ on $2 < \alpha < 3$ is equivalent to the constraint 
$\frac{3}{4} (17^{\frac{1}{2}}-1) < \alpha < 3$. This indicates that the condition $u \in L^{\infty}(L^{\alpha ,\infty})$ with such a $\alpha$ lying in 
$ ( \frac{3}{4} (17^{\frac{1}{2}}-1) , 3)$ is beyond the classical Prodi-Serrin-Ladyzhenskaya range and the $L^{\infty}(L^{3})$ criterion of \cite{Esca}. 
This means that this extra assumption, which is the first tehnical condtion in the hypothesis of Theorem \ref{goal}, is reasonable.\

Next, let us mention that there is nothing deep about the second condition in Theorem \ref{goal}, which says that
the large velocity region $\{x \in \mathbb{R}^{3} : |u(t, x )| > M_{0}\}$ of the solution is restricted within the open ball $\{x \in \mathbb{R}^{3} : |x| < r_{0}\}$
for some given radius $r_{0}$. Even without this second condition as in the hypothesis of Theorem \ref{goal}, the same qualitive result which says that large velocity
region of a solution $u$ to \eqref{NS} has to be within a certain ball with some sufficiently large radius $R_{0}$ depending on $u$ can be deduced by means of an application
of the partial regularity theorem of Caffarelli, Kohn and Nirenberg. So, in this sense, the second condition imposed on Theorem \ref{goal} is not very crucial and is 
imposed only for convenience.\

Before we finish the discussion about Theorem \ref{goal}, we point out that the proof of Theorem \ref{goal} as presented in section 4 of our paper closely follows the 
proof of the regularity criterion in \cite{smoothnesscriteria}. However, we also point out that we have given completely \emph{new} idea which allows us to use the extra
weak $L^{\alpha}$ space condition with $\frac{3}{4} (17^{\frac{1}{2}} - 1) < \alpha < 3$ to lower down the index $\beta$ in \eqref{nonlinearrecurrence} 
from the typical value of $\frac{5}{3}$ or 
$\frac{4}{3}$ to become as close to $1$ as possible. For those readers who are interested only in those new ideas contributed to the proof of Theorem \ref{goal}, 
we have given, in Section 3 of our paper, an outline of those crucial and important ideas which make the old argument of \cite{smoothnesscriteria} become 
strong enough to arrive at Theorem \ref{goal}. But we also give, in Section 4 of our paper, the complete details of the proof of Theorem \ref{goal} by including those
new ideas outlined in Section 3 in the technical argument.    

Besides the main result in Theorem \ref{goal}, we also will construct, in section 2, a 
divergence free velocity field $u$ within a stream-tube segment 
with increasing twisting (i.e. increasing swirl) among the streamlines of $u$ towards the ending cross section of the stream-tube. 
The construction of such a velocity field $u$ as in section 2 
demonstrates the way in which \emph{excessive} twisting of streamlines towards the ending cross section of the stream-tube can result in the \emph{blow up} of
the quantities $\|u\|_{L^{\alpha}(\mathbb{R}^{3})}$ (for some $2 < \alpha < 3$) and $\|\dv (\frac{u}{|u|})\|_{L^{6}(\mathbb{R}^{3})}$ while at the same time preserving
the finite energy property $u \in L^{2}(\mathbb{R}^{3})$ of the fluid. We \emph{do not} claim that this explicit construction of such a divergence free vector field has 
anything to do with actual solutions
to \eqref{NS}. One purpose of such a construction of a divergence free velocity field $u$ with increasing swirl towards the ending cross section of the stream-tube is to 
illustrate the possibility of having a finite energy velocity field with increasing swirl which is beyond the scope covered by the regularity criterion of Vasseur in
\cite{Vasseur2} and the $L^{\infty}(L^{3})$ criterion of \cite{Esca}. 
In a certain sense, the excessive twisting of streamlines of the velocity field as constructed in Section 2 within a stream-tube
segment with \emph{almost constant cross section everywhere} (see Definition \ref{uniformlybundle}) will cause the streamlines to become \emph{densely packed} together 
towards the ending cross section of the stream-tube, and this denser and denser packing of streamlines eventually leads to the blow up of the velocity field at a singular 
point lying at the center of the ending cross section of the stream-tube. According to the regularity criterion in Theorem \ref{goal}, one can speculate that if 
the velocity field $u$ as constructed in Section 2 with increasing swirl towards the ending cross section of the stream-tube can be \emph{realized} as 
\emph{an instantaneous profile $v(T,\cdot )$} of a time-dependent solution $v : [0,T)\times \mathbb{R}^{3} \rightarrow \mathbb{R}^{3}$ to \eqref{NS} in which 
singularity occurs at the blow up time $T$, then, it must be that 
the rate of increase of $F = \dv (\frac{u}{|u|})$ along those streamlines with increasing twisting must go \emph{beyond} the expontential growth rate. 
Even though the construction in Section 2 is very interesting, it is totally independent of the regularity criterion of Theorem \ref{goal}, and 
the reader should treat this as a separate topic.  
    
\section{Possible blow-up situation of a velocity field with large swirl}
In this section, we \emph{attempt} to characterize a divergence free velocity field $u$ which is 
specified in a stream-tube segment around a representative streamline (with an incoming cross section and an ending cross section) in a such a way that the streamlines 
generated by such a velocity field will have unbounded increasing swirl 
(ie increasing twisting around the representative streamline) \emph{towards the ending cross section of the stream-tube segment}.
The uncontrolled increasing swirl of streamlines towards the ending cross section of the stream-tube segment will lead to an isolated singularity
located at the point of intersection between the center representative streamline and the ending cross section of the stream-tube.
Based on the above consideration, we will give \emph{necessary} conditions which characterize the properties $u \in L^{2}$, $u\not \in  L^\alpha$ (for some given $2<\alpha<3$)
and $\text{div}\frac{u}{|u|}\not\in  L^6$ of such a velocity field with increasing swirl.\\ 

In order to describe such a velocity field $u$ with increasing swirl towards the ending cross section of the stream-tube segment, we first specify the center representative
streamline $\gamma_{\eta} : [0,S)\rightarrow \mathbb{R}^{3}$ as follow.
\begin{defn}
 (Representative stream line.)

Let $\gamma_\eta : [0,S) \rightarrow \mathbb{R}^{3}$ be such that
\begin{equation}
  \partial_s\gamma_\eta(s)=\frac{u}{|u|}(\gamma_\eta(s))\quad\text{and}\quad\gamma_\eta(0)=\eta
  \in \mathbb{R}^{3}.
\end{equation}
\end{defn}
Note that the ending value $S$ is excluded from the definition of the representative streamline, since $\gamma_{\eta}(S)$ is supposed to be the isolated singularity point 
created by the unbounded increasing swirl of those streamlines close to the representative streamline. Before we can create the stream-tube segment around the representative
streamline $\gamma_{\eta}$, we need to specify the initial streamplane $A$ with parameter r as follow.
\begin{defn}
 (Initial stream plane $A$ with parameter $r$.)
Let $\{\bar A_0(r)\}_{r\in (0,1]}$ be a smooth family of smoothly bounded open set in $\mathbb{R}^2$ s.t.
$\bar A_0(r)\subset \bar A_0(r')$ ($r<r'$),  $\bar A_0(r)\to \{0\}$ ($r\to 0$).
Let 
\begin{equation}
 A_0(r)=\{x\in\mathbb{R}^3: R(x-\eta)\in \bar A_0(r)\},
\end{equation}
where $R$ is a rotation matrix s.t. $R(\frac{u}{|u|}(\eta,t))=(0,0,1)$.
\end{defn}
The initial streamplane $A_{0}(1)$ is exactly the incoming cross section of the stream-tube segment which will be specified. To construct the stream-tube segment
with $A_{0}(1)$ as its incoming cross section, we just need to specify, for each $s \in (0,S)$, the associated stream-plane $A(r,s)$ intersecting $\gamma_{\eta}$ at 
the point $\gamma_\eta(s)$ as follow. 
\begin{defn}\label{Streamplanes}
(Stream-planes.)
Let
\begin{multline*}
 A(r,s):=\bigcup_{\eta'\in A_0(r)}\{\gamma_{\eta'}(s'):
 s'\  \text{is the minimum among all possible}\ \tau>0
 \ \text{for which}\ \\
  \gamma_{\eta'}(\tau)
 \text{belongs the plane which passes through the point}\  \gamma_\eta(s) \ \text{and is perpendicular to}\  
 \partial_s\gamma_\eta(s) \}
\end{multline*}
\end{defn}
For simplicity, we just set $A(s):=A(1,s)$. Then, we can define the stream-tube to be 
\begin{equation}
T_{[0,S)}^{A} = \bigcup_{0\leqslant s < S}A(s) 
\end{equation}
We remark that, for any $x\in A(s)$, there is $r$ s.t. $x\in\partial A(r,s)$. This is due to the fact that for each $s \in [0,S)$, $A(r,s)$ is strictly shrinking towards the 
representative streamline as $r\rightarrow 0^{+}$. Based on this observation, we introduce an orthonormal coordinate frame system within the stream-tube $T_{[0,S)}^{A}$ in the
following definition.
\begin{defn}
 For $x,y\in\partial A(r,s)$, let $e_\theta(x):=\lim_{y\to x}\frac{x-y}{|x-y|}$, $e_z(x):=\frac{u}{|u|}(\gamma_\eta(s))$ and 
let $e_r(x)$ be s.t. 
\begin{equation}
 \langle e_z(x),e_r(x)\rangle=\ \langle e_\theta(x),e_r(x)\rangle=0\quad\text{and}\quad |e_r(x)|=1.
\end{equation}
\end{defn}
We emphasize that the notations $e_\theta$, $e_z$ and $e_r$ are borrowed from the notations of the standard cylindrical coordinate frame $\partial_{r}$, 
$\frac{1}{r}\partial_{\theta}$ and $\partial_{z}$ for axisymmetric velocity field about the $z$-axis. This is a good choice of notation, since one can imagine 
that the representative streamline $\gamma_\eta$ plays a similar role as the axi-symmetric axis provided $\gamma_\eta$ is \emph{relatively straight}.\
Next, in order to describe the increasing swirl of $u$ towards the ending cross section $A(S) = A(1,S)$ of the stream-tube segment $T_{[0,S)}^{A}$, we will now decompose 
$\frac{u}{|u|}$ into its radial component, $z$-component, and swirl component as in the following definition. 
\begin{defn}
 (Decomposition of normalized streamline.)

Let $\omega_\theta$, $\omega_r$ and $\omega_z$ be s.t.
\begin{equation}
 \frac{u}{|u|}(x)=\omega_\theta(x)e_\theta(x)+\omega_r(x)e_r(x)+\omega_z(x)e_z(x).
\end{equation}
\end{defn}

\begin{remark}\label{omega_z}
We see that $\omega_\theta^2+ \omega_r^2+ \omega_z^2=1$ and  $\omega_z(x)\to 1$ ($x\to\gamma_\eta(s)$) if $u$ is smooth.
\end{remark}
In order to give a model of possible blow-up situation,  we need to define ``uniformly bundle" as follows:

\begin{defn}\label{uniformlybundle}
We call that ``the stream-tube segment $T_{[0,S)}^{A}$ has a uniformly bundle" if the following two properties hold
\begin{itemize}
\item For any $B(0)\subset A(0)$ and any $s\in[0,S]$, we have $C^{-1}\leqslant\frac{|B(s)|}{|B(0)|}\leqslant C$, for some universal constant $C > 0$.
Here, $B(s)$ is defined in the same way as $A(r,s)$ through replacing $A_{0}(r)$ by $B(0)$ in Definition \ref{Streamplanes}.
\item For the same universal constant $C> 0$, we have $\sup_{y\in A(0)}u\cdot e_z(y)\leq C\inf _{y\in A(0)}u\cdot e_z(y)$.
\end{itemize}


\end{defn}

\begin{remark}
Since $\int_{B(0)}u\cdot e_z(y)d\sigma_y=\int_{B(s)}u\cdot e_z(y)d\sigma_y$ by divergence
free, we see $u\cdot e_z(x)=\lim_{B(s)\ni x}\frac{1}{|B(s)|}\int_{B(s)}u\cdot e_z(y)d\sigma_y
\approx\lim_{B(0)\ni x'}\frac{1}{|B(0)|}\int_{B(0)}u\cdot e_z(y)d\sigma_y=u\cdot e_z(x')$
for any two points $x\in A(s)$ and $x'\in A(0)$ connected by a streamline passing through $A(0)$ and $A(s)$ , if $A(s)$ has a uniformly bundle.
\end{remark}
\begin{remark}
  If $A(s)$ has a uniformly bundle, we can see from  divergence free condition
\begin{equation}
\int_{B(0)}u\cdot e_z(y)d\sigma_y=\int_{B(s)}u\cdot e_z(y)d\sigma_y\approx |B(s)||u\cdot e_z|=|B(s)||u|\omega_z
\end{equation}
and then
\begin{equation}\label{calfromdivfree}
 \frac{\int_{B(0)}u\cdot e_z(y)d\sigma_y}{|B(s)||\omega_z(x)|}\approx|u(x)|
 \quad\text{for}\quad x\in B(s)\subset A(s).
\end{equation}
\end{remark}

Now, we want to characterize the properties $u\in L^2$, $u\not \in  L^\alpha$ (for some given $2<\alpha<3$) and $\text{div}\frac{u}{|u|}\not\in  L^6$ in terms of 
some conditions specifying \emph{how fast the streamlines are increasing their swirl towards the ending corss section $A(s)$ of the stream-tube segment $T_{[0,S)}^{A}$}.\\


To specify the increasing swirl of streamlines towards the ending cross section $A(S)$ of the stream-tube segment $T_{[0,S)}^{A}$, we need to decompose each stream-plane $A(s)$ 
into the disjoint union of a countable list of ring-shaped regions $A_j(s)$ as follow. We first select a decreasing sequence of positive numbers 
$\{r_j\}_{j= 1}^{\infty}$ dropping down to $0$($r_j\searrow 0$) as $j \rightarrow \infty$. We then set $A_j(s):=A(r_j,s)\setminus A(r_{j+1},s)$. Notice that $A_{j}(s)$ 
is shrinking towards the representative streamline $\gamma_{\eta}$ as $j$ becomes large. We also set
 \begin{equation}
 \omega_z^{A_j}(s):=\int_{A_j(s)}\omega_z(y)d\sigma_y/|A_j(s)|.
 \end{equation}
That is, $\omega_{z}^{A_{j}}(s)$ is the average of $\omega_{z}$ over the ring-shaped region $A_{j}(s)$ in the stream-plane $A(s)$. We can assume, according to Definition
\ref{uniformlybundle}, that 
\begin{equation}\label{roughlysamecross}
 |A_j(s)|\approx|A_j(0)|\quad\text{for}\quad s\in (0,S].
\end{equation}
Since we require that $u$ blows up at the isolated singluar point $\gamma_{\eta}(S)$ lying in $A(S)$, in light of condition \eqref{roughlysamecross} and 
\eqref{calfromdivfree}, we would require that, as $s$ becomes close to $S$, $\omega_{z}^{A_j}(s)$ should become small as $j$ becomes large, which indicates that the average
swirl (or twisting) of those streamlines passing through $A_{j}(s)$ should become large as $s\rightarrow S$ and $j\rightarrow \infty$.

Now, in order to ensure that $u\in  L^2$, we impose $(S-s)^{\frac{1}{2}-\epsilon}$ as the lower bound for $\omega_{z}$ as follow.
\begin{itemize}
\item (The condition to ensure $u\in  L^2$ ) For any $0 < s < S$, we have $(S-s)^{\frac{1}{2}-\epsilon} < \omega_{z}(\gamma_{\eta'}(s)) < 1$ for any $\eta'\in A(0)$.
\end{itemize}

The purpose of the above condition is to prevent the  swirl of  streamlines passing through $A_{j}(s)$ to become \emph{too large} as $s\rightarrow S$ and 
$j\rightarrow \infty$, because we want to have the finite energy property for $u$.
Under the condition $(S-s)^{\frac{1}{2}-\epsilon} < \omega_{z}(\gamma_{\eta'}(s)) < 1$ and \eqref{calfromdivfree}, a direct calculation yields the finite $L^{2}$ property of $u$ as follow.
\begin{equation}
\|u\|_{L^2}^2\approx \|u\|_{L^2(T^A_{[0,S)})}^2=\int_{T^A_{[0,S)}}\left|\frac{\omega^A(0)}{\omega_z(x)}\right|^2dx
\approx\int_0^S\int_{A}\left|\frac{\omega^A(0)}{\omega_z(\gamma_\eta(s))}\right |^2d\eta ds
 \leqslant \int_0^S\frac{C}{(S-s)^{1-2\epsilon}}ds<\infty.
\end{equation}
 In order to ensure that $u\not\in  L^\alpha$, we impose $(S-s)^{\frac{1}{\alpha}}$ as the upper bound for $\omega_{z}$ as follow.
 \begin{itemize}
 \item (The condition to ensure $u\not\in  L^\alpha$ )
Let $\{S_j\}_j\subset [0,S)$ be s.t. $S_j\to S$ ($j\to\infty$). 
For any $0<s<S_j$, we have 
\begin{equation*}
|\omega_z(\gamma_{\eta'}(s))|\leq  (S-s)^{1/\alpha}
\quad\text{for}\quad \eta'\in A_j(0).
\end{equation*}
\end{itemize}
We show $u\not\in L^\alpha$.  By Remark \ref{omega_z}, we see $\omega^{A_j}(0)\approx 1$. Thus 
\begin{equation}
\|u\|_{L^\alpha}^\alpha\approx \|u\|_{L^\alpha(T^A_{[0,S)})}^\alpha=\int_{T^A_{[0,S)}}\left|\frac{\omega^A(0)}{\omega_z(x)}\right|^\alpha dx
\approx
\sum_j\int_0^S\int_{A_j}\left|\frac{\omega^{A_j}(0)}{\omega_z(\gamma_\eta(s))}\right |^\alpha d\eta ds
 \geqslant \int_0^S\frac{C}{(S-s)}ds=\infty.
\end{equation}

In order to show $\left\|\text{div}\frac{u}{|u|}\right\|_{L^6} = \infty$, 
 we impose $(S-s)^{-1}$ as the upper bound for $|\partial_s\omega^{A_j}_{z}(s)|^6$ as follow.\begin{itemize}
\item (The condition to ensure $\text{div}(u/|u|)\not\in  L^6$ )
Let $\{\tilde S_j\}_j\subset [0,S)$ be s.t. $S_j<\tilde S_j <S$ and 
\begin{equation*}
|A_j(0)|\int_{S_j}^{\tilde S_j}(S-s)^{-1}ds\geqslant C,
\end{equation*}
where $C$ is a universal constant.
For any $S_j<s<\tilde S_j$, we have 
\begin{equation*}
|\partial_s\omega^{A_j}_z(s)|^6>(S-s)^{-1}.
\end{equation*}
\end{itemize}
We need to get a \emph{rough} expression of 
$\frac{1}{A_{j}(s)}\int_{A_{j}(s)}\text{div}(\frac{u}{|u|}) dy$
as follow. Let $s > 0$ be fixed. Then, for any $s_{1} > s$ to be sufficiently close to $s$, we consider the following stream-tube $T_{[s,s_{1}]}^{A_{j}}$ connecting the 
stream-plane $A_j(s)$ to $A_j(s_{1})$. 
\begin{equation}
T_{[s,s_{1}]}^{A_j} = \bigcup_{ s \leqslant \tau \leqslant s_{1} } A_j(s). 
\end{equation}
From Definition \ref{Streamplanes}, we can view the stream-tube $T_{[s,s_{1}]}^{A_j}$ as being formed by the union of those streamlines which first pass into the 
stream-tube through the cross section $A_j(s)$ and eventually leave the same stream-tube through the cross section $A_j(s_{1})$. Since $s_{1}$ is chosen to be close to $s$, the
stream-tube $T_{[s,s_{1}]}^{A_j}$ is roughly the same as the product $A_j(s)\times [s , s_{1}]$, which, together with condition \eqref{roughlysamecross}, makes the following 
deduction justifiable. 
\begin{eqnarray*}
 \frac{1}{A_{j}(s)}\int_{A_j(s)}\text{div}(\frac{u}{|u|}) dy&=&\lim_{s_{1}\rightarrow s}\frac{1}{(s_{1} - s)}\int_{s}^{s_{1}} 
\frac{1}{|A_j(\tau )|}\int_{A_j(\tau)}\text{div}(\frac{u}{|u|})dy d\tau \\
&\approx & \lim_{s_{1}\rightarrow s} \frac{1}{(s_{1} - s)|A_j(s)|} \int_{T_{[s,s_{1}]}^{A_j}} \text{div}(\frac{u}{|u|})dy\\
&=&\lim_{s_{1} \rightarrow s}\frac{1}{|A_j(s)|(s_{1} - s)} \{\int_{A_j(s_1)}\frac{u}{|u|}\cdot e_zd\sigma- \int_{A_j(s)}\frac{u}{|u|}\cdot e_zd\sigma\}\\
&=&\lim_{s_{1} \rightarrow s} \frac{1}{|A_j(s)|(s_{1} - s)} \{\int_{A_j(s_{1})} \omega_zd\sigma-\int_{A_j(s)}\omega_zd\sigma \}\\
&=&\lim_{s_{1} \rightarrow s}\frac{1}{|A_j(s)|(s_1-s)}\left(\omega^{A_j}_z(s_1)|A_j(s_1)|-\omega_z^{A_j}(s)|A_j(s)|\right)\\
&=&\frac{1}{|A_j(s)|}\partial_s\{\omega_z^{A_j}(s)|A_j(s)|\} =\left(\partial_s\omega_z^{A_j}(s)+\frac{\partial_s|A_j(s)|}{|A_j(s)|}\omega_z^{A_j}(s)\right).
\end{eqnarray*}
Hence, it follows from the above calculation and an application of Holder inequality that

\begin{equation*}
\begin{split}
\int_{A_{j}(s)} |\dv(\frac{u}{|u|})|^{6} &\geqslant \frac{1}{|A_{j}(s)|^{5}} \left| \int_{A_{j}(s)} \dv(\frac{u}{|u|}) \right|^{6}\\ 
& \approx |A_{j}(s)| \left(\partial_s\omega_z^{A_j}(s)+\frac{\partial_s|A_j(s)|}{|A_j(s)|}\omega_z^{A_j}(s)\right)^{6} \\
& \approx |A_{j}(0)| \left(\partial_s\omega_z^{A_j}(s)+\frac{\partial_s|A_j(s)|}{|A_j(s)|}\omega_z^{A_j}(s)\right)^{6}\\
&\geq |A_j(0)||\partial_s\omega_z^{A_j}(s)|^6. 
\end{split}
\end{equation*}
Therefore, 

\begin{equation*}
\begin{split}
\left\|\dv\frac{u}{|u|}\right\|_{L^6(T^A_{[0,S)})}^6\approx
\sum_j\int_{[0,S)}\int_{A_{j}(s)} |\dv(\frac{u}{|u|})|^{6} 
&\geq \sum_j\int_{[S_j,\tilde S_j)}|A_j(0)||\partial_s\omega_z^{A_j}(s)|^6ds\geqslant \sum_j C=\infty.\\
\end{split}
\end{equation*}



\section{Outline of the proof of Theorem \ref{goal}}

The proof of Theorem \ref{goal} is quite similar to the one in \cite{smoothnesscriteria}. The purpose of this section is just to outline those \emph{crucial and 
important} changes which have to be made to the structure of the proof as pesented in \cite{smoothnesscriteria}, so that the modified proof will be strong enough
to give the result of Theorem \ref{goal}. In other words, we will only state the essential changes to the main argument in \cite{smoothnesscriteria}
 which are the new ideas contributed in this paper.
  
Just in the same way as \cite{smoothnesscriteria}, we will follow the parabolic De Giorgi's method devolped by Vasseur in \cite{Vasseur}. 
So, let us fix our notation as follow. We remark that, without the lost of generality, we will assume that the possible blow up time $T$ is just $1$.\ 

\begin{itemize}
\item for each $k \geqslant 0$, let $Q_{k} = [T_{k} , 1]\times \mathbb{R}^3$, in which
$T_{k} = \frac{3}{4} - \frac{1}{4^{k+1}}$.

\item for each $k \geqslant 0$, let $v_{k} = \{ |u| -R(1 - \frac{1}{2^k}) \}_{+}$.
\item for each $k \geqslant 0$, let $w_{k} = \{ |u| -R^{\beta}(1 - \frac{1}{2^k}) \}_{+}$ , with $\beta > 1$ to be selected later.
\item for each $k \geqslant 0$, let $d_{k}^{2} = \frac{R( 1 - \frac{1}{2^k})}{|u|}
\chi_{\{|u| > R (1- \frac{1}{2^k})\}} |\nabla |u||^{2} + \frac{v_{k}}{|u|} |\nabla
u|^2$.
\item for each $k \geqslant 0$, let $D_{k}^{2} = \frac{R^{\beta}( 1 - \frac{1}{2^k})}{|u|}
\chi_{\{|u| > R^{\beta} (1- \frac{1}{2^k})\}} |\nabla |u||^{2} + \frac{w_{k}}{|u|} |\nabla
u|^2$.
\item for each $k \geqslant 0$, let $U_{k} = \frac{1}{2}\|v_{k}\|^{2}_{L^{\infty}(T_{k} , 1 ;
L^{2}(\mathbb{R}^3))} + \int_{T_{k}}^{1}\int_{\mathbb{R}^{3}} d_{k}^2 dx\,dt$.
\end{itemize}
With the above setting, the first author proved the following proposition (see \cite{smoothnesscriteria}).
\begin{prop}\label{index}
Let $u$ be a suitable weak solution for the Navier-Stokes equation on 
$[0,1]\times \mathbb{R}^3$ which satisfies the condition that
$|\frac{u\cdot \nabla F}{|u|^{\gamma }}|\leqslant A |F|$, where $A$ is some finite-positive constant, and $\gamma $ is some positive number satisfying 
$0 < \gamma < \frac{1}{3}$. Then, there exists some constant $C_{p,\beta }$, depending only on $1<p<\frac{5}{4}$, and $\beta > \frac{6-3p}{10-8p}$,and also some constants $0 < \alpha , K <\infty$, which do depend on our suitable weak solution $u$, such that the following inequality holds

\begin{equation}\label{conclusionprop2.1}
\begin{split}
U_{k}\leqslant 
& C_{p,\beta }2^{\frac{10k}{3}}
\{\frac{1}{R^{\beta \frac{10-8p}{3p}-\frac{2-p}{p}}} \|u\|_{L^{\infty}(0,1;L^{2}(\mathbb{R}^{3}))}^{2(1-\frac{1}{p})}U_{k-1}^{\frac{5-p}{3p}} +\\
&(1+A)(1+\frac{1}{\alpha })
(1+K^{1-\frac{1}{p}})(1+\|u\|_{L^{\infty}(0,1;L^{2}(\mathbb{R}^{3}))})\times \\
&[ (\frac{1}{R^{\frac{10}{3}-2p\beta +1-\gamma -p }})^{\frac{1}{p}}U_{k-1}^{\frac{5}{3p}}
 + \frac{1}{R^{\frac{10}{3}-2\beta -\gamma }}U_{k-1}^{\frac{5}{3}}]\} ,
\end{split}
\end{equation}
for every sufficiently large $R>1$.
\end{prop}

The nonlinear recurrence relation as given in \eqref{conclusionprop2.1} was indeed the main cornerstone leading to the regularity criterion in \cite{smoothnesscriteria}.
More precisely, the structure of \eqref{conclusionprop2.1} directly gives the smallness of $U_{1}$ as long as $R$ is sufficiently large. The smallness of $U_{1}$,
together with the nonlinear recurrence structure of relation \eqref{conclusionprop2.1}, then 
allowed us to deduce in \cite{smoothnesscriteria} the decay of $U_{k}$ to $0$ (as $k \rightarrow \infty$) by means of the following useful lemma as appeared in 
\cite{Vasseur}.
\begin{lemma}\label{Vass}
For any given constants $B$, $\beta > 1$, there exists some constant
$C^*_{0}$ such that for any sequence $\{a_{k}\}_{k\geqslant 1}$
satisfying $0 < a_{1} \leq C^*_{0}$ and $a_{k} \leqslant B^k
a^{\beta}_{k-1}$, for any $k \geqslant 1$, we have
$lim_{k\rightarrow \infty} a_{k} = 0$ .
\end{lemma}
The resulting decay of $U_{k}$ to $0$ as $k \rightarrow \infty$ allowed the first author to draw the conclusion that $u$ is essentially bounded by some sufficiently
large constant $R > 1$ over $[\frac{3}{4} , 1) \times \mathbb{R}^{3}$, and this lead to the following theorem in \cite{smoothnesscriteria}.

\begin{thm}\label{goal0}
Let $u : [0,T)\times\mathbb{R}^{3} \rightarrow \mathbb{R}^{3}$ be a Leray-Hopf solution to \eqref{NS} 
which is smooth on $[0,T)\times\mathbb{R}^{3}$ (with $T$ to be the possible blow up time) and which satisfies the condition that 
$|\frac{u\cdot \nabla F }{|u|^{\gamma }}|\leqslant A|F|$, in which $A$ is some positive constant, and $\gamma $ is some positive constant for which
$0< \gamma < \frac{1}{3}$. Then, it follows that the $u$ is $L^{\infty}$-bounded on $[\frac{3}{4} , 1)\times \mathbb{R}^{3}$ and hence the
smoothness of $u$ can be extended beyond $T$.
\end{thm}

In this paper we will refine the $\gamma$ in Theorem \ref{goal0} to be $1$.
As indicated in the introduction, the problem we face here is that those powers of $U_{k-1}$ such as 
$\frac{5-p}{3p}$, $\frac{5}{3p}$ and $\frac{5}{3}$ (appearing in Proposition \ref{index}), are \emph{too far} from $1$.
However, the use of Lemma \ref{Vass} only requires that $\beta>1$, so the extra condition 
$u \in L^{\infty}(0,1;L^{\alpha ,\infty}(\mathbb{R}^{3}))$, with $\alpha \in (2,3)$ satisfying $1+ 2(\frac{\alpha}{3} + \frac{3}{\alpha}) > 0$ can help us 
to bring the powers of $U_{k-1}$ to become very close to $1$, and this in turn allows us to replace the old condition $|\frac{u\cdot \nabla F}{|u|^{\gamma}}| \leq A |F|$  
with $\gamma \in (0,\frac{1}{3})$ by the new one $|\frac{u\cdot \nabla F}{|u|}| \leq A |F|$.



Technically speaking, the key idea which allows us to use the condition 
$u \in L^{\infty}(L^{\alpha ,\infty})$ (with $\alpha \in (2,3)$ satisfying $1+ 2(\frac{\alpha}{3} + \frac{3}{\alpha}) > 0$) to lower down the powers
of $U_{k-1}$ to become close to $1$ is the following lemma. 
We can establish the following lemma for any truncations $w_{k-1} = (|u| - R^{\beta}(1-\frac{1}{2^{k-1}}))_{+}$ (with $k \geq 2$) of a Leray-Hopf 
solution $u \in L^{\infty}(0, 1 ; L^2(\mathbb{R}^3)) \cap L^2(0, 1 ; \dot{H}^1(\mathbb{R}^3))$ satsifying 
the condition $u \in L^{\infty}(0, 1 ; L^{\alpha , *}(\mathbb{R}^3))$ for some given $\alpha \in (2, 3)$.

\begin{lemma}\label{weakLP}
Consider a Leray-Hopf weak solution $u \in L^{\infty}(0, 1 ; L^2(\mathbb{R}^3)) \cap L^2(0, 1 ; \dot{H}^1(\mathbb{R}^3))$ which satsifies the condition $u \in L^{\infty}(0, 1 ; L^{\alpha , *}(\mathbb{R}^3))$ for some given $\alpha \in (2, 3)$. Then, the truncation $w_{k-1} = (|u| - R^{\beta}(1-\frac{1}{2^{k-1}}))_{+}$ of $|u|$ satisfies the following inequality for each $k \geqslant 2$ and each $\delta$ with $0 < \delta < \frac{4}{3}$. 

\begin{equation}\label{weakLPimprove1}
\int_{Q_{k-1}} w_{k-1}^{\frac{10}{3}} \leqslant C_{0} \{ \frac{2^{\alpha -1}}{\alpha - 2} \|u\|_{L^{\infty}(0, 1 ; L^{\alpha , *}(\mathbb{R}^3))}\}^{\frac{2}{3} - \delta}
\frac{U_{k-1}^{1+ \delta}}{R^{\beta (\alpha - 2) (\frac{2}{3} - \delta)}} ,
\end{equation}

in which $C_{0}$ is a universal constant essentially arising from the Sobolev embedding theorem. In the same way, the truncation $v_{k} = (|u| - R(1-\frac{1}{2^{k}}))_{+}$ also satisfies the following inequality for each $k \geqslant 2$ and each $\delta$ with $0 < \delta < \frac{4}{3}$. 

\begin{equation}\label{weakLPimprove2}
\int_{Q_{k-1}} v_{k-1}^{\frac{10}{3}} \leqslant C_{0} \{ \frac{2^{\alpha -1}}{\alpha - 2} \|u\|_{L^{\infty}(0, 1 ; L^{\alpha , *}(\mathbb{R}^3))}\}^{\frac{2}{3} - \delta}
\frac{U_{k-1}^{1+ \delta}}{R^{ (\alpha - 2) (\frac{2}{3} - \delta)}} .
\end{equation}

\end{lemma}

\begin{proof}

To begin, let $u \in L^{\infty}(0, 1 ; L^2(\mathbb{R}^3)) \cap L^2(0, 1 ; \dot{H}^1(\mathbb{R}^3))$ to be a Leray-Hopf solution which satsifies the 
condition $u \in L^{\infty}(0, 1 ; L^{\alpha , *}(\mathbb{R}^3))$ for some given $\alpha$ with $2< \alpha < 3 $. 
Recall that the truncation $w_{k-1} = (|u| - R^{\beta}(1-\frac{1}{2^{k-1}}))_{+}$
satisfies the property that $|\nabla w_{k-1}| \leqslant D_{k-1} \leqslant 5^{\frac{1}{2}} d_{k-1}$, for $k \geqslant 2$ (The relation 
$|\nabla w_{k-1}| \leqslant D_{k-1}$ can be verified easily, while the relation $D_{k-1} \leqslant 5^{\frac{1}{2}} d_{k-1}$ was justified in Lemma 4.1 of \cite{smoothnesscriteria}). 
So, it follows from standard interpolation inequality that

\begin{equation}
\begin{split}
\int_{Q_{k-1}} w_{k-1}^{\frac{10}{3}} &\leqslant C_{0} \|w_{k-1}\|_{L^{\infty}(T_{k-1} , 1 ; L^{2}(\mathbb{R}^{3}))}^{\frac{4}{3}} \|\nabla w_{k-1}\|_{L^{2}(Q_{k-1})}^{2}\\
& \leqslant C_{0} \{ sup_{t \in [T_{k-1} , 1]} \int_{\mathbb{R}^{3}} w_{k-1}^{2}(t,x) dx      \}^{\frac{2}{3}} U_{k-1}\\
& \leqslant C_{0} U_{k-1}^{1+ \delta} \{ sup_{t \in [T_{k-1} , 1]} \int_{\mathbb{R}^{3}} w_{k-1}^{2}(t,x) dx      \}^{\frac{2}{3} - \delta} .
\end{split}
\end{equation}

But according to the assumption that $u \in L^{\infty}(0, 1 ; L^{\alpha , *}(\mathbb{R}^3))$, we can control $\int_{\mathbb{R}^{3}} w_{k-1}^{2}(t,x) dx$ (for each $k \geqslant 2$) uniformlly over $t \in [0,1]$ as follow.

\begin{equation}
\begin{split}
\int_{\mathbb{R}^{3}} w_{k-1}^{2}(t,x) dx & = 2 \int_{0}^{\infty} r |\{x \in \mathbb{R}^{3} : w_{k-1}(t,x) > r\}| dr \\
& \leqslant 2 \int_{0}^{\infty} r |\{ x \in \mathbb{R}^{3} : |u(t,x)| > r + R^{\beta}(1- \frac{1}{2^{k-1}})     \}| dr \\
& \leqslant 2 \int_{0}^{\infty} (r + \frac{R^{\beta}}{2}) |\{ x \in \mathbb{R}^{3} : |u(t,x)| > r +  \frac{R^{\beta}}{2}     \}| dr \\
& = 2 \int_{\frac{R^{\beta}}{2}}^{\infty} r  |\{ x \in \mathbb{R}^{3} : |u(t,x)| > r      \}| dr\\
& \leq 2  \|u\|_{ L^{\infty}(0, 1 ; L^{\alpha , *}(\mathbb{R}^3))} \int_{\frac{R^{\beta}}{2}}^{\infty} r^{1- \alpha } dr \\
& = \frac{2^{\alpha -1}}{\alpha -2} \|u\|_{ L^{\infty}(0, 1 ; L^{\alpha , *}(\mathbb{R}^3))} \frac{1}{R^{\beta (\alpha -2)}} .
\end{split}
\end{equation}

Hence, inequality (\ref{weakLPimprove1}) follows from the above two inequality estimations. By the same way, we can also derive inequality (\ref{weakLPimprove2}) by replacing $w_{k}$ by $v_{k} = (|u| - R(1-\frac{1}{2^{k}}))$ and $R^{\beta}$ by $R$.

\end{proof}

As a corollary of Lemma~\ref{weakLP}, we have the following result which allows us to raise up the index for the terms  $\|\chi_{\{ w_{k} > 0 \} }\|_{L^q(Q_{k-1})}$
and $\|\chi_{\{ v_{k} > 0 \} }\|_{L^q(Q_{k-1})}$.

\begin{lemma}\label{cheb}
Suppose that the given suitable weak solution $u : [0,1]\times \mathbb{R}^{3} \rightarrow \mathbb{R}$ satisfies the condition $u \in L^{\infty}(0, 1 ; L^{\alpha ,*}(\mathbb{
R}^{3}))$ for some given $\alpha \in (2,3)$. Then, for any $1 < q < \infty$, and any $1 < \delta <  \frac{4}{3}$, we have

\begin{equation}\label{improvedchebineqone}
\|\chi_{\{ w_{k} > 0 \} }\|_{L^q(Q_{k-1})} \leqslant C_{(\alpha , \delta , q)} \frac{2^{\frac{10k}{3q}}}{R^{\frac{1}{q}[\frac{10\beta}{3} + 
\beta (\alpha - 2)(\frac{2}{3} - \delta)] }} \cdot \|u\|_{L^{\infty}(L^{\alpha , *})}^{(\frac{2}{3}-\delta )\frac{1}{q}} U_{k-1}^{(1+ \delta)\frac{1}{q}} ,
\end{equation} 
in which the constant $C_{(\alpha , \delta , q)}$ is given by $ C_{(\alpha , \delta , q)} = C_{0}^{\frac{1}{q}} [\frac{2^{\alpha -1}}{(\alpha -2)}]^{(\frac{2}{3}- \delta)\frac{1}{q}} $, 
with $C_{0}$ to be a universal constant arising from the Sobolev embedding theorem and standard interpolation.\
  
In the same way, we have the following estimate for $\|\chi_{\{ v_{k} > 0 \} }\|_{L^q(Q_{k-1})}$, with $1< q < \infty$ and $1 < \delta < \frac{4}{3} $.
\begin{equation}\label{improvedchebineqtwo}
\|\chi_{\{ v_{k} > 0 \} }\|_{L^q(Q_{k-1})} \leqslant C_{(\alpha , \delta , q)} \frac{2^{\frac{10k}{3q}}}{R^{\frac{1}{q}[\frac{10}{3} + 
 (\alpha - 2)(\frac{2}{3} - \delta)] }} \cdot \|u\|_{L^{\infty}(L^{\alpha , *})}^{(\frac{2}{3}-\delta )\frac{1}{q}} U_{k-1}^{(1+ \delta)\frac{1}{q}} .
\end{equation} 

\end{lemma}

\emph{Remark} Notice that the constant  $C_{(\alpha , \delta , q)}$ as appears in inequality (\ref{improvedchebineqone}) blows up to $\infty$ as the choice 
of $\alpha$ approaches to $2$, which means that inequality (\ref{improvedchebineqone}) applies only in the case of $\alpha > 2$. We also point out that replacing the 
old Lemma 3.2 and Lemma 3.3 in \cite{smoothnesscriteria} by the above lemma (i.e. Lemma \ref{cheb}) is the \emph{crucial} decision leading to the
final success of our new proof of Theorem \ref{goal} (see the next section, in which we will give all the details of the new proof of Theorem \ref{goal} ).  
\begin{proof}
We recall that the sequence of truncations $w_{k}$ is defined to be $w_{k} = (|u| - R^{\beta}(1-\frac{1}{2^{k}}))_{+}$. So, it is easy to see that $\{w_{k} > 0\} \subset \{w_{k-1} > \frac{R^{\beta}}{2^{k}}\}$. Hence, it follows from inequality (\ref{weakLPimprove1}) that 

\begin{equation}
\begin{split}
\int_{Q_{k-1}} \chi_{\{ w_{k} > 0  \}} & \leqslant   \int_{Q_{k-1}} \chi_{\{ w_{k-1} > \frac{R^{\beta}}{2^{k}}  \}} \\
& \leqslant \frac{2^{\frac{10k}{3}}}{R^{\frac{10\beta}{3}}} \int_{Q_{k-1}} w_{k-1}^{\frac{10}{3}} \\
& \leqslant \frac{2^{\frac{10k}{3}}}{R^{\frac{10\beta}{3}}} \cdot C_{0} \{ \frac{2^{\alpha -1}}{\alpha - 2} \|u\|_{L^{\infty}(0, 1 ; L^{\alpha , *}(\mathbb{R}^3))}\}^{\frac{2}{3} - \delta} \frac{U_{k-1}^{1+ \delta}}{R^{\beta (\alpha - 2) (\frac{2}{3} - \delta)}} .
\end{split}
\end{equation}

Hence, inequality (\ref{improvedchebineqone}) follows from taking the power $\frac{1}{q}$ on both sides of the above inequlity. The deduction of inequality 
\eqref{improvedchebineqtwo} follows in the same way.
\end{proof}

In order to adopt to the new hypothesis $|u\cdot \nabla F| \leqslant A|u|\cdot|F|$ on $\{(t,x) \in [0,1)\times B(r_{0}) : |F(t,x)| \geqslant L\}$ (for some given constant $L> 0$)
, the second refinement is on the function $\psi$ appearing in Step five of the proof in \cite{smoothnesscriteria}.
We \emph{redefine} the function $\psi : \mathbb{R}\rightarrow \mathbb{R}$ as the one which satisfies the following conditions 
\begin{itemize}
\item $\psi (t)=1$, for all $t\geqslant L+1$.
\item $0 < \psi (t) < 1$, for all $t$ with $L< t < L+1$.
\item $\psi (t) = 0$, for all $-L \leqslant t \leqslant L$ .
\item $-1 < \psi (t) < 0$, for all $t$ with $-L-1 < t < -L$.
\item $\psi (t) = -1$, for all $t \leqslant -L-1$.
\item $0 \leqslant \frac{d}{dt}\psi  \leqslant 2$, for all $t\in \mathbb{R}$. 
\end{itemize}

We further remark that the smooth function $\psi : \mathbb{R}\rightarrow \mathbb{R}$ characterized by the above properties must also satisfy the property that
$ \frac{d \psi }{dt}|_{(t)} = 0  $, on 
$t \in (-\infty , -L-1)\cup(-L,L)\cup \cup (L+1 , \infty )$.

Up to this point, we have already spelled out \emph{all} the important changes that have to be made to the old argument in \cite{smoothnesscriteria}. In the next section, 
we will redo the old argument in \cite{smoothnesscriteria} by including all those important changes given here, and see the way in which the modified new argument 
will lead to the result of Theorem \ref{goal}.

\section{Appendix: Technical steps of the proof of Theorem \ref{goal}.}

\vskip0.2cm

The purpose of this section is to convince the readers of the correctness of the outline in the previous section through giving all the technical details of the proof of 
Theorem \ref{goal}. Except those \emph{crucial and important} changes as given in the outline of the previous section, 
the structure of the proof of theorem \ref{goal} is in many aspects the same as the one in \cite{smoothnesscriteria}. It is also not suprising that some of the 
technical aspects of the proof of Theorem \ref{goal} as given below are directly transported (or copied) from that of \cite{smoothnesscriteria} (This is justified for those parts to which
no change is necessary). So, in a certain sense, all the new ideas of the proof of Theorem \ref{goal} has already been given in the outline of the previous section, 
and we spell out all the
details of the proof of Theorem \ref{goal} here only for the sake of completeness. Moreover, we remark that, within this section,
 the definitions of $T_{k}$, $Q_{k}$, $v_{k}$, $w_{k}$ $d_{k}$ etc were given in the beginning of section 3. Moreover, the possible finite blow up time for the solution
$u ; [0,1)\times \mathbb{R}^{3} \rightarrow \mathbb{R}^{3}$ under consideration is assumed to be $1$.\\

\noindent{\bf Step one}

To begin the argument, we recall that, according to Lemma 5 in \cite{Vasseur}, the truncations $v_{k} = \{|u| - R (1-\frac{1}{2^{k}})\}$ of a given suitable weak 
solution $u : [0,1]\times \mathbb{R}^{3} \rightarrow \mathbb{R}^{3}$ satisfy the following inequality in the sense of distribution. 
\begin{equation}\label{LEIfortruncations}
\partial_{t} (\frac{v_{k}^2}{2}) + d_{k}^{2} - \triangle (\frac{v_{k}^2}{2}) + \dv(\frac{v_{k}^2}{2} u) +
\frac{v_{k}}{|u|} u \nabla P \leqslant 0 .
\end{equation}
Next, let us consider the variables $\sigma$ , $t$ verifying $T_{k-1} \leqslant \sigma \leqslant T_{k} \leqslant t
\leqslant 1$. Then, we have
\begin{itemize}
\item $\int_{\sigma}^{t} \int_{\mathbb{R}^3} \partial_{t} (\frac{v_{k}^2}{2}) dx\,ds =
\int_{\mathbb{R}^3} \frac{v_{k}^2(t,x)}{2} dx - \int_{\mathbb{R}^3} \frac{v_{k}^2(\sigma ,x)}{2} dx$.
\item $\int_{\sigma}^{t}\int_{\mathbb{R}^3} \triangle (\frac{v_{k}^2}{2}) dx\,ds = 0$.
\item $\int_{\sigma}^{t} \int_{\mathbb{R}^3} div (\frac{v_{k}^2}{2} u) dx\,ds = 0$.
\end{itemize}
So, it is straightforward to see that
\begin{equation*}
\int_{\mathbb{R}^3} \frac{v_{k}^2(t,x)}{2} dx + \int_{\sigma}^{t}\int_{\mathbb{R}^3}d_{k}^2dx\,ds
\leqslant \int_{\mathbb{R}^3} \frac{v_{k}^2(\sigma ,x)}{2}dx +
\int_{\sigma}^{t}    \vert \int_{\mathbb{R}^3}\frac{v_{k}}{|u|} u \nabla P dx \vert       ds ,
\end{equation*}
for any $\sigma$, $t$ satisfying $T_{k-1}\leqslant \sigma \leqslant T_{k} \leqslant t \leqslant 1$. By taking the
average over the variable $\sigma$, we yield
\begin{equation*}
\int_{\mathbb{R}^3} \frac{v_{k}^2(t,x)}{2} dx + \int_{T_{k}}^{t}\int_{\mathbb{R}^3}d_{k}^2 dx\,ds
\leqslant \frac{4^{k+1}}{6}   \int_{T_{k-1}}^{T_{k}}\int_{\mathbb{R}^3}v_{k}^2(s,x)dx\,ds
+ \int_{T_{k-1}}^{t}|\int_{\mathbb{R}^3}\frac{v_{k}}{|u|} u \nabla Pdx|ds .
\end{equation*}
By taking  the sup over $t\in [T_{k} , 1]$. the above inequality will give the following
\begin{equation*}
U_{k} \leqslant \frac{4^{k+1}}{6}   \int_{Q_{k-1}}v_{k}^2 + \int_{T_{k-1}}^{1}|\int_{\mathbb{R}^3}\frac{v_{k}}{|u|}u\nabla
Pdx|ds .
\end{equation*}
But, by using the interpolation inequality $\|f\|_{L^{\frac{10}{3}}(Q_{k})} \leqslant \|f\|_{L^{\infty}(T_{k},1; L^{2}(\mathbb{R}^{3}))}^{\frac{2}{5}}
\|\nabla f\|_{L^{2}(Q_{k})}^{\frac{3}{5}}$ (see Lemma 3.1 of \cite{smoothnesscriteria} or \cite{Vasseur}  ) and the inequality 
$\|\chi_{v_{k} > 0}\|_{L^{q}(Q_{k-1})} \leqslant (\frac{2^{k}}{R})^{\frac{10}{3q}} C^{\frac{1}{q}} U_{k-1}^{\frac{5}{3q}}$ (see Lemma 3.2 of \cite{smoothnesscriteria}
or \cite{Vasseur}), we can carry out the following estimate.

\begin{equation*}
\begin{split}
\int_{Q_{k-1}}v_{k}^2 & = \int_{Q_{k-1}}v_{k}^2 \chi_{\{v_{k} > 0\}} \\
& \leqslant (\int_{Q_{k-1}}v_{k}^{\frac{10}{3}})^{\frac{3}{5}} \|\chi_{\{v_{k} > 0 \}}\|_{L^{\frac{5}{2}}(Q_{k-1})}\\
& \leqslant \|v_{k}\|_{L^{\frac{10}{3}}(Q_{k-1})}^{2} \frac{2^{\frac{4k}{3}}}{R^{\frac{4}{3}}}      C^{\frac{2}{5}}U_{k-1}^{\frac{2}{3}}\\
& \leqslant \|v_{k-1}\|_{L^{\frac{10}{3}}(Q_{k-1})}^{2} \frac{ 2^{\frac{4k}{3}}}{R^{\frac{4}{3}}}           
C^{\frac{2}{5}}U_{k-1}^{\frac{2}{3}}\\
& \leqslant CU_{k-1}^{\frac{5}{3}} \frac{ 2^{\frac{4k}{3}}}{R^{\frac{4}{3}}} .
\end{split}
\end{equation*}
As a result, we have the following conclusion
\begin{equation}\label{trival}
U_{k}\leqslant \frac{ 2^{\frac{10k}{3}}}{R^{\frac{4}{3}}}     C U_{k-1}^{\frac{5}{3}} +
\int_{T_{k-1}}^{1}|\int_{\mathbb{R}^3}\frac{v_{k}}{|u|}u \nabla p dx|ds .
\end{equation}

\noindent{\bf Step two}

Now, in order to estimate the term $\int_{T_{k-1}}^{1}|\int_{\mathbb{R}^3}\frac{v_{k}}{|u|}u\nabla Pdx|ds$, we would
like to carry out the following computation
\begin{equation*}
\begin{split}
-\triangle P &= \sum \partial_{i} \partial_{j} (u_{i}u_{j})\\
& = \sum \partial_{i} \partial_{j} \{ (1-\frac{w_{k}}{|u|})u_{i} (1-\frac{w_{k}}{|u|})u_{j} \}
+  2\sum \partial_{i} \partial_{j} \{ (1-\frac{w_{k}}{|u|}) u_{i} \frac{w_{k}}{|u|} u_{j}  \}\\
& + \sum \partial_{i} \partial_{j} \{ \frac{w_{k}}{|u|} u_{i} \frac{w_{k}}{|u|} u_{j} \} ,
\end{split}
\end{equation*}
in which $w_{k}$ is given by $w_{k} = \{|u| - R^{\beta} (1- \frac{1}{2^{k}} )\}_{+}$, and $\beta > 1$ is some arbritary index which will be determined later. This 
motivates us to decompose $P$ as $P = P_{k1} + P_{k2} + P_{k3} $, in which
\begin{equation}\label{P1}
-\triangle P_{k1} = \sum \partial_{i} \partial_{j} \{ (1-\frac{w_{k}}{|u|})u_{i}
 (1-\frac{w_{k}}{|u|})u_{j}  \} ,
\end{equation}

\begin{equation}\label{P2}
- \triangle P_{k2} = \sum \partial_{i} \partial_{j} \{ 2(1-\frac{w_{k}}{|u|})u_{i} \frac{w_{k}}{|u|}
 u_{j} \}
\end{equation}

\begin{equation}\label{P3}
-\triangle P_{k3} = \sum \partial_{i} \partial_{j}\{\frac{w_{k}}{|u|}u_{i} \frac{w_{k}}{|u|}u_{j} \} .
\end{equation}

Here, we have to remind ourself that the cutting functions which are used in the decomposition of the pressure are indeed 
$w_{k} = \{|u| - R^{\beta} (1- \frac{1}{2^{k}}) \}_{+}$, for all $k \geqslant 0$ , in which $\beta$ is some suitable index strictly greater than $1$. With respect 
to the cutting functions $w_{k}$, we need to define the respective $D_{k}$ as follow:
\begin{equation*}
D_{k}^{2} = \frac{R^{\beta} (1-\frac{1}{2^{k}})}{|u|} \chi_{\{ w_{k} > 0 \}} |\nabla |u||^{2}
+ \frac{w_{k}}{|u|} |\nabla u|^{2} .
\end{equation*}
Then, just like what happens to the cutting functions $v_{k}$, we have the following assertions about the cutting functions $w_{k}$, which are easily verified 
(see \cite{Vasseur}).
\begin{itemize}
\item $|\nabla w_{k}| \leqslant D_{k}$, for all $k \geqslant 0$.
\item $|\nabla (\frac{w_{k}}{|u|}u_{i})| \leqslant 3 D_{k}$, for all $k \geqslant 0$, and $1 \leqslant i \leqslant 3$.
\item $ |\nabla (\frac{w_{k}}{|u|}) u_{i} | \leqslant 2 D_{k} $, for any $k \geqslant 0$, and
$1 \leqslant i \leqslant 3 $.
\item $D_{k} \leq 5^{\frac{1}{2}} d_{k}$ as long as $R$ is larger than some fixed constant $R_{0}$ (see Lemma 4.1 of \cite{smoothnesscriteria} for a proof of this).
\end{itemize}




Now, let us recall that we have already used the cutting functions $w_{k}$ to obtain the decomposition $P = P_{k1} + P_{k2} + P_{k3}$, in which 
$P_{k1}$, $P_{k2}$, and$P_{k3}$ are described in equations (\ref{P1}), (\ref{P2}), and (\ref{P3}) respectively.\\
Due to the incompressible condition $\dv (u) = 0$, we have the following two identities
\begin{itemize}
\item $ \int_{\mathbb{R}^{3}}  \frac{v_{k}}{|u|}  u \nabla P_{k2} dx    = \int_{\mathbb{R}^{3}}  ( \frac{v_{k}}{|u|} - 1  ) u \nabla P_{k2} dx  $.
\item $ \int_{\mathbb{R}^{3}}  \frac{v_{k}}{|u|}  u \nabla P_{k3} dx    = \int_{\mathbb{R}^{3}}  ( \frac{v_{k}}{|u|} - 1  ) u \nabla P_{k3} dx  $.
\end{itemize}
Hence, it follows that
\begin{equation}\label{Divide}
\begin{split}
\int_{T_{k-1}}^{1} |\int_{\mathbb{R}^{3}} \frac{v_{k}}{|u|} u \nabla P dx  |dt
&\leqslant \int_{T_{k-1}}^{1} |\int_{\mathbb{R}^{3}} \nabla (\frac{v_{k}}{|u|}) u P_{k1} dx |dt +\int_{Q_{k-1}} (1-\frac{v_{k}}{|u|})|u| |\nabla P_{k2}|\\
& + \int_{Q_{k-1}} (1-\frac{v_{k}}{|u|})|u| |\nabla P_{k3}| .
\end{split}
\end{equation}

\noindent{\bf Step 3}

We are now ready to deal with the term  $\int_{Q_{k-1}}(1-\frac{v_{k}}{|u|})|u||\nabla P_{k2}|$. For this purpose, let $p$ be such that $1< p < \frac{5}{4} $, and let 
$q = \frac{p}{p-1}$, so that $2<q<\infty$. We remark that the purpose of the condition $1< p < \frac{5}{4} $ is to ensure that the quantity $\frac{2p}{2-p}$ will satisfy the 
condition $2 <  \frac{2p}{2-p} < \frac{10}{3}$, which is required in the forthcoming inequality estimation \ref{6required}. Next, by applying Holder's inequality, we find 
that
\begin{equation*}
\begin{split}
\|(1-\frac{v_{k}}{|u|})u \|_{L^{q}(\mathbb{R}^{3})}
& \leqslant
\|(1-\frac{v_{k}}{|u|})u\|_{L^{2}(\mathbb{R}^{3})}^{\frac{2}{q}}
\|(1-\frac{v_{k}}{|u|})u\|_{L^{\infty}(\mathbb{R}^{3})}^{1-\frac{2}{q}}\\
&\leqslant R^{1-\frac{2}{q}} \|(1-\frac{v_{k}}{|u|})u\|_{L^{2}(\mathbb{R}^{3})}^{\frac{2}{q}}\\
&\leqslant  R^{\frac{2}{p}-1} 
\|u\|_{L^{\infty}(0,1;L^{2}(\mathbb{R}^{3}))}^{2(1-\frac{1}{p})}
\end{split}
\end{equation*}
Hence, it follows from Holder's inequality that
\begin{equation*}
\int_{\mathbb{R}^{3}} (1-\frac{v_{k}}{|u|})|u||\nabla P_{k2}| dx 
\leqslant R^{\frac{2}{p}-1} 
\|u\|_{L^{\infty}(0,1;L^{2}(\mathbb{R}^{3}))}^{2(1-\frac{1}{p})} 
\{\int_{\mathbb{R}^{3}} |\nabla P_{k2}|^{p} dx \}^{\frac{1}{p}} .
\end{equation*}
Hence, we have
\begin{equation}\label{4}
\int_{Q_{k-1}} (1-\frac{v_{k}}{|u|})|u||\nabla P_{k2}| 
\leqslant R^{\frac{2}{p}-1}
\|u\|_{L^{\infty}(0,1;L^{2}(\mathbb{R}^{3}))}^{2(1-\frac{1}{p})}
\|\nabla P_{k2}\|_{L^{p}(Q_{k-1})} .
\end{equation}
But, we recognize that 
\begin{equation*}
\nabla P_{k2} = \sum R_{i}R_{j}\{2(1-\frac{w_{k}}{|u|})u_{i} \nabla [\frac{w_{k}}{|u|}u_{j}] 
+ 2(1-\frac{w_{k}}{|u|})u_{j}[\frac{w_{k}}{|u|}\nabla u_{i}]
- 2 \nabla [\frac{w_{k}}{|u|}]u_{i} \frac{w_{k}}{|u|}u_{j}\}.
\end{equation*}
Moreover, it is straightforward to see that for any $1 \leqslant i,j \leqslant 3$, we have
\begin{itemize}
\item $ |2(1-\frac{w_{k}}{|u|})u_{i} \nabla [\frac{w_{k}}{|u|}u_{j}] +  2(1-\frac{w_{k}}{|u|})u_{j}[\frac{w_{k}}{|u|}\nabla u_{i}]| \leqslant 8 R^{\beta } D_{k}$.
\item $|2 \nabla [\frac{w_{k}}{|u|}]u_{i} \frac{w_{k}}{|u|}u_{j}| \leqslant 8 w_{k} D_{k}$.
\end{itemize}
So, we can decompose $\nabla P_{k2}$ as $\nabla P_{k2} = G_{k21} + G_{k22}$, where $G_{k21}$ and $G_{k22}$ are given by
\begin{itemize}
\item $G_{k21} = \sum R_{i} R_{j} \{2(1-\frac{w_{k}}{|u|})u_{i} \nabla [\frac{w_{k}}{|u|}u_{j}] +   2(1-\frac{w_{k}}{|u|})u_{j}[\frac{w_{k}}{|u|}\nabla u_{i}]\}$.
\item $G_{k22} = -\sum R_{i} R_{j} \{2 \nabla [\frac{w_{k}}{|u|}]u_{i} \frac{w_{k}}{|u|}u_{j}\}$.
\end{itemize}
In order to use inequality (\ref{4}), we need to estimate $\|G_{k21}\|_{L^{p}(Q_{k-1})}$ and $\|G_{k22}\|_{L^{p}(Q_{k-1})}$ respectively, for $p$ with $1 < p < \frac{5}{4} $. 
Indeed, by applying the Zygmund-Calderon Theorem, we can deduce that
\begin{itemize}
\item $\|G_{k21}\|_{L^{p}(Q_{k-1})} \leqslant C_{p} R^{\beta } \|D_{k}\|_{L^{p}(Q_{k-1})}$,
\item $ \|G_{k22}\|_{L^{p}(Q_{k-1})} \leqslant C_{p} \|w_{k}D_{k}\|_{L^{p}(Q_{k-1})}$,
\end{itemize}
where $C_{p}$ is some constant depending only on $p$. But it turns out that
\begin{equation*}
\begin{split}
\|D_{k}\|_{L^{p}(Q_{k-1})}^{p} & = \int_{Q_{k-1}} D_{k}^{p} \chi_{\{w_{k} > 0 \}}\\
& \leqslant \{\int_{Q_{k-1}}D_{k}^{2}\}^{\frac{p}{2}} 
\|\chi_{\{w_{k} > 0 \}}\|_{L^{\frac{2}{2-p}}(Q_{k-1})}\\
&\leqslant   5^{\frac{p}{2}}   \|d_{k}\|_{L^{2}(Q_{k-1})}^{p}
C_{\alpha , p} \frac{2^{\frac{5(2-p)k}{3}}}{R^{\beta (\frac{2-p}{2})[\frac{10}{3} + (\alpha - 2)(\frac{2}{3} - \delta)] }} \cdot 
\|u\|_{L^{\infty}(L^{\alpha , *})}^{(\frac{2}{3}-\delta )(\frac{2-p}{2})} U_{k-1}^{(1+ \delta)(\frac{2-p}{2})} \\
& \leqslant C_{\alpha , p ,\delta } \frac{2^{\frac{5(2-p)k}{3}}}{R^{\beta (\frac{2-p}{2})[\frac{10}{3} +  (\alpha - 2)(\frac{2}{3} - \delta)] }} \cdot 
\|u\|_{L^{\infty}(L^{\alpha , *})}^{(\frac{2}{3}-\delta )(\frac{2-p}{2})} \cdot U_{k-1}^{1+ \delta (\frac{2-p}{2})} .
\end{split}
\end{equation*}
That is , we have
\begin{equation*}
\|D_{k}\|_{L^{p}(Q_{k-1})} \leqslant  C_{\alpha , p ,\delta } \frac{2^{\frac{5(2-p)k}{3p}}}{R^{\beta (\frac{2-p}{2p})[\frac{10}{3} + (\alpha - 2)(\frac{2}{3} - \delta)] }} 
\cdot \|u\|_{L^{\infty}(L^{\alpha , *})}^{(\frac{2}{3}-\delta )(\frac{2-p}{2p})} \cdot U_{k-1}^{\frac{1}{p} + \delta (\frac{2-p}{2p})} .
\end{equation*}
Hence, it follows that
\begin{equation}\label{5}
\|G_{k21}\|_{L^{p}(Q_{k-1})} \leqslant     C_{\alpha , p ,\delta } \frac{2^{\frac{5(2-p)k}{3p}}}{R^{\beta [\frac{10-8 p}{3p} + (\frac{2-p}{2p})(\alpha - 2)(\frac{2}{3} - 
\delta )] }} \cdot \|u\|_{L^{\infty}(L^{\alpha , *})}^{(\frac{2}{3}-\delta )(\frac{2-p}{2p})} \cdot U_{k-1}^{\frac{1}{p} + \delta (\frac{2-p}{2p})}  .
\end{equation}
On the other hand, we have
\begin{equation*}
\begin{split}
\|w_{k}D_{k}\|_{L^{p}(Q_{k-1})}^{p} & = \int_{Q_{k-1}} w_{k}^{p} D_{k}^{p}\\
& \leqslant \{ \int_{Q_{k-1}}w_{k}^{\frac{2p}{2-p}}\}^{\frac{2-p}{2}}
\{\int_{Q_{k-1}}D_{k}^2  \}^{\frac{p}{2}}\\
& \leqslant C_{p} \{\int_{Q_{k-1}} w_{k}^{\frac{2p}{2-p}}\}^{\frac{2-p}{2}} 
U_{k-1}^{\frac{p}{2}} .
\end{split}
\end{equation*}
Now, let us recall that $1 < p < \frac{5}{4}$, and put $r = \frac{2p}{2-p}$. we then recognize that $2 < r = \frac{2p}{2-p} < \frac{10}{3}$, if $1 < p < \frac{5}{4}$. So, we 
can have the following estimation
\begin{equation}\label{6required}
\begin{split}
\int_{Q_{k-1}}w_{k}^{\frac{2p}{2-p}}
 & = \int_{Q_{k-1}} w_{k}^{r} \chi_{\{ w_{k} > 0  \}}\\
&\leqslant \int_{Q_{k-1}} w_{k}^{r} \chi_{\{ w_{k-1} > \frac{R^{\beta}}{2^{k}}\}}\\
&\leqslant \frac{1}{R^{\beta (\frac{10}{3} -r)}} 2^{k( \frac{10}{3} - r )}
\int_{Q_{k-1}} w_{k}^{\frac{10}{3}}\\
&\leqslant \frac{C_{\alpha , \delta} \|u\|_{L^{\infty}(L^{\alpha , *})}^{\frac{2}{3} - \delta}}{R^{\beta [\frac{20-16p}{3(2-p)} + (\alpha - 2)(\frac{2}{3} - \delta )  ]}} 
2^{\frac{k(20-16p)}{3(2-p)}} 
U_{k-1}^{1+\delta} .
\end{split}
\end{equation}
Hence, it follows that
\begin{equation}\label{6}
\begin{split}
\|G_{k22}\|_{L^{p}(Q_{k-1})} & \leqslant C_{p} \|w_{k}D_{k}\|_{L^{p}(Q_{k-1})} \\
& \leqslant C_{\alpha, p , \delta }
\frac{2^{ \frac{(10-8p)k}{3p}} }{R^{\beta [\frac{10-8p}{3p} + 
(\frac{2-p}{2p})(\alpha - 2) (\frac{2}{3} -\delta )]}} \|u\|_{L^{\infty}(L^{\alpha , *})}^{ (\frac{2}{3}-\delta )(\frac{2-p}{2p})} U_{k-1}^{\frac{1}{p} + 
\delta (\frac{2-p}{2p})} .
\end{split}
\end{equation}
By combining inequalities (\ref{4}), (\ref{5}), (\ref{6}), we deduce that
\begin{equation}\label{middle}
\begin{split}
\int_{Q_{k-1}} (1-\frac{v_{k}}{|u|})|u||\nabla P_{k2}|
\leqslant  \frac{2^{ \frac{(10-8p)k}{3p}} C( \alpha , p, \delta ; u )  }{R^{\beta [\frac{10-8p}{3p} + (\frac{2-p}{2p})(\alpha - 2) (\frac{2}{3} -\delta )] 
-(\frac{2-p}{p})}}  U_{k-1}^{\frac{1}{p} + \delta (\frac{2-p}{2p})} ,
\end{split}
\end{equation}
in which the constant $C( \alpha , p, \delta ; u )$ is in the form of 
\begin{equation}\label{specialconstant}
C( \alpha , p, \delta ; u ) = C_{\alpha , p , \delta } \|u\|_{L^{\infty}(L^{2})}^{2(1-\frac{1}{p})} 
\|u\|_{L^{\infty}(L^{\alpha , *})}^{ (\frac{2}{3}-\delta )(\frac{2-p}{2p})} .
\end{equation}
As for the term $\int_{Q_{k-1}} (1-\frac{v_{k}}{|u|})|u||\nabla P_{k3}|$. We first notice that
\begin{equation*}
P_{k3} = \sum R_{i}R_{j}\{\frac{w_{k}}{|u|}u_{i}\frac{w_{k}}{|u|}u_{j}\} . 
\end{equation*}
So, we know that 
\begin{equation*}
\nabla P_{k3} = \sum R_{i}R_{j} \{\nabla [\frac{w_{k}}{|u|}u_{i}]\frac{w_{k}}{|u|}u_{j}
+ \frac{w_{k}}{|u|}u_{i} \nabla [\frac{w_{k}}{|u|} u_{j} ] \} ,
\end{equation*}
with 
\begin{equation*}
| \nabla [\frac{w_{k}}{|u|}u_{i}]\frac{w_{k}}{|u|}u_{j}
+ \frac{w_{k}}{|u|}u_{i} \nabla [\frac{w_{k}}{|u|} u_{j} ] |
\leqslant 6 w_{k}D_{k} .
\end{equation*}
Again, by the Risez's theorem, we have  $\|\nabla P_{k3}\|_{L^{p}(\mathbb{R}^{3})} \leqslant C_{p} 
\|w_{k}D_{k}\|_{L^{p}(\mathbb{R}^{3})}$, in which $C_{p}$ is some constant depending only on $p$. So, we can repeat the same type of estimation, just as what we have done to 
the term $\int_{Q_{k-1}}(1-\frac{v_{k}}{|u|})|u||\nabla P_{k2}|$, to conclude that

\begin{equation}\label{easy}
\begin{split}
\int_{Q_{k-1}} (1-\frac{v_{k}}{|u|})|u||\nabla P_{k3}|
&\leqslant R^{\frac{2}{p}-1}\|u\|_{L^{\infty}(0,1;L^{2}(\mathbb{R}^{3}))}^{2(1-\frac{1}{p})}
\|\nabla P_{k3}\|_{L^{p}(Q_{k-1})}\\
&\leqslant \frac{2^{ \frac{(10-8p)k}{3p}} C( \alpha , p, \delta ; u )  }{R^{\beta [\frac{10-8p}{3p} + (\frac{2-p}{2p})(\alpha - 2) (\frac{2}{3} -\delta )] -(\frac{2-p}{p})}} 
 U_{k-1}^{\frac{1}{p} + \delta (\frac{2-p}{2p})} ,
\end{split}
\end{equation}
in which the constant $C( \alpha , p, \delta ; u ) $ is again in the form of (\ref{specialconstant}).\\

We have to ensure that the quantity $\beta [\frac{10-8p}{3p} + (\frac{2-p}{2p})(\alpha - 2) (\frac{2}{3} -\delta )] -(\frac{2-p}{p})   $ is \emph{strictly greater than $0$}. 
To this end, recall that $p > 1 $ can be as close to $1$ as possible, and $\delta > 0$ can also be as close to $0$ as possible. So, by passing to the limit as 
$p\rightarrow 1^{+}$, and $\delta \rightarrow 0^{+}$, we have
\begin{equation}\label{limitingvalueofindexofR}
\lim_{p\rightarrow 1^{+} ,\delta \rightarrow 0^{+} } \beta [\frac{10-8p}{3p} + (\frac{2-p}{2p})(\alpha - 2) (\frac{2}{3} -\delta )] -(\frac{2-p}{p}) = 
\beta(\frac{\alpha}{3}) - 1.
\end{equation}
Now, we insist that the choice of $\beta$ has to satisify the condition $\beta > \frac{3}{\alpha}$, under which we must have the limiting value $\beta(\frac{\alpha}{3}) - 1$ 
to be strictly positive. Hence, for such a choice of $\beta$, it follows from (\ref{limitingvalueofindexofR}) that the following relation holds for all $p > 1$ to be 
sufficiently close to $1$, and all $\delta > 0$ to be sufficiently close to $0$.
\begin{equation}
 \beta [\frac{10-8p}{3p} + (\frac{2-p}{2p})(\alpha - 2) (\frac{2}{3} -\delta )] -(\frac{2-p}{p}) > 0.
\end{equation}

\noindent{\bf Step four}

We now have to raise up the index for the term $\int_{T_{k-1}}^{1}|\int_{\mathbb{R}^{3}} \nabla (\frac{v_{k}}{|u|})uP_{k1}dx|ds$.\

Recall that, in the hypothesis of Theorem \ref{goal}, there is some constant $M_{0} > 0$ for which $|u| \leqslant M_{0}$ is valid on the outer region 
$[0,1)\times \{x \in \mathbb{R}^{3} : |x| \geqslant r_{0}\}$ for some given radius $r_{0} > 0$. As a result, we will now choose $R > 2M_{0}$ so that, for each 
$k\geq 1$ and $t\in [0,1)$ we have $\{|u(t,\cdot )| > R(1-\frac{1}{2^{k}})\} \subset B(r_{0})$, which means that both $v_{k}(t, \cdot )$ and $d_{k}(t,\cdot )$ are compactly 
supported in $B(r_{0})$. Hence, for such a choice of $R > 2M_{0}$, we always can express $U_{k}$ as
\begin{equation*}
U_{k} = \frac{1}{2} \sup_{t\in [T_{k}, 1)} \int_{B(r_{0})} v_{k}^{2}(t,\cdot ) dx + \int_{T_{k}}^{1} \int_{B(r_{0})} d_{k}^{2} dx dt .
\end{equation*} 
Since $\nabla (\frac{v_{k}}{|u|})u = -R(1-\frac{1}{2^{k}})F\chi_{\{v_{k} > 0\}} $, we have \emph{for any} $R > 2M_{0} $ that
\begin{equation*}
\begin{split}
|\int_{\mathbb{R}^{3}}\nabla (\frac{v_{k}}{|u|})u P_{k1}dx|
&= |\int_{B(r_{0})} R(1-\frac{1}{2^{k}})F\chi_{\{v_{k} > 0 \}}P_{k1}dx|\\
&\leqslant R\int_{B(r_{0})}|F|\chi_{\{v_{k} > 0 \}} |P_{k1} - (P_{k1})_{B(r_{0})}|dx\\
&+  R\int_{B(r_{0})}|F|\chi_{\{v_{k} > 0\}}|(P_{k1})_{B(r_{0})}|dx ,
\end{split} 
\end{equation*}
for all $k\geqslant 1$, and all $\frac{1}{2} < t < 1$ (here, the symbol $(P_{k1})_{B}$ stands for the average value of
 $P_{k1}$ over the ball $B$ ).
From now on, we will \emph{always assume, within this section, that our choice of $R$ has to satisify $R > 2M_{0}$ }. Now, since 
$P_{k1} = \sum R_{i}R_{j} \{(1-\frac{w_{k}}{|u|})u_{i}(1-\frac{w_{k}}{|u|})u_{j}\}$, it follows from the Risez's Theorem in the theory of singular integral that
$\|P_{k1}(t,\cdot )\|_{L^{2}(\mathbb{R}^{3})} 
\leqslant C_{2}R^{\beta }\|u(t,\cdot )\|_{L^{2}(\mathbb{R}^{3})}$, for all $t\in [0,1]$, in which $C_{2}$ is some constant depending only on $2$. So, we can use the 
Holder's inequality to carry out the following estimation
\begin{equation*}
\begin{split}
|(P_{k1})_{B(r_{0})}(t)| & \leqslant\frac{1}{|B(r_{0})|}\int_{B(r_{0})}|P_{k1}(t,x)|dx\\ 
&\leqslant \frac{1}{|B(r_{0})|^{\frac{1}{2}}} \|P_{k1}(t,\cdot )\|_{L^{2}(B(r_{0}))}\\
&\leqslant \frac{1}{|B(r_{0})|^{\frac{1}{2}}}C_{2}R^{\beta }\|u(t,\cdot )\|_{L^{2}(\mathbb{R}^{3})}\\
& \leqslant C(r_{0}) R^{\beta }\|u\|_{L^{\infty}(0,1; L^{2}(\mathbb{R}^{3}))} ,
\end{split}
\end{equation*}
in which the constant $C(r_{0}) = \frac{1}{|B(r_{0})|^{\frac{1}{2}}}C_{2}$ depends on $r_{0}$. As a result, it follows that
\begin{equation}\label{10}
\begin{split}
|\int_{\mathbb{R}^{3}}\nabla (\frac{v_{k}}{|u|}) u P_{k1}dx|
& \leqslant R \int_{B(r_{0})}|F|\chi_{\{v_{k} > 0\}}|P_{k1}- (P_{k1})_{B(r_{0})}|dx\\
& +  C(r_{0})R\|u\|_{L^{\infty}(0,1;L^{2}(\mathbb{R}^{3}))} 
\int_{B(r_{0})}R^{\beta } |F|\chi_{\{v_{k} > 0\}}
\end{split}
\end{equation}
Indeed, the operator $R_{i}R_{j}$ is indeed a Zygmund- Calderon operator, and so $R_{i}R_{j}$ must be a bounded operator from 
$L^{\infty}(\mathbb{R}^{3})$ to $BMO(\mathbb{R}^{3})$. Hence we can deduce that
\begin{equation*}
\begin{split}
\|P_{k1}(t, \cdot )- (P_{k1})_{B(r_{0})}(t)\|_{BMO}
& = \|P_{k1}(t,\cdot )\|_{BMO}\\
&\leqslant C_{0} 
\|(1-\frac{w_{k}}{|u|})u_{i}(1-\frac{w_{k}}{|u|})u_{j}\|_{L^{\infty}(\mathbb{R}^{3})}\\
&\leqslant C_{0}R^{2\beta } ,
\end{split} 
\end{equation*}
for all $t\in (0,1)$, in which $C_{0}$ is some constant depending only on $\mathbb{R}^{3}$.\\

Just as the proof of the main result in \cite{smoothnesscriteria}, at this stage, we need the assistant of the following Lemma, which is a straightforward corollary of the famous
$BMO$ result \cite{John} of John and Nirenberg. For a proof of this lemma, we refer to Lemma 4.3 of \cite{smoothnesscriteria}.

\begin{lemma}\label{BMO}
(see \cite{smoothnesscriteria})Let $B$ be a ball with finite radius sitting in $\mathbb{R}^{3}$. There exists some finite positive constants $\alpha$ and $K$,depending only on $B$, 
such that for every $\mu \geqslant 0$, every $f\in BMO(\mathbb{R}^3)$ with $\int_{B}fdx =0$, and $p$ with $1 < p < \infty$, we have
$\int_{B} \mu |f| \leqslant \frac{2p}{\alpha (p-1)} 
\{ 1+ K^{1-\frac{1}{p}}\}  \|f\|_{BMO}
\{ (\int_{B}\mu )^{\frac{1}{p}}  +  \int_{B} \mu log^{+}\mu  \}$.
\end{lemma}
   
So, we now apply Lemma~\ref{BMO} with $\mu = |F|\chi_{\{v_{k} > 0 \}}$, and 
$f = P_{k1} - (P_{k1})_{B(r_{0})}$ to deduce that
\begin{equation*}
\begin{split}
\int_{B(r_{0})} |F| \chi_{\{v_{k} > 0 \}}|P_{k1}-(P_{k1})_{B(r_{0})}|dx
&\leqslant \frac{2pC_{0}}{\alpha (p-1)} \{1+ K^{1-\frac{1}{p}}\}\times \\
&\{(\int_{B(r_{0})}R^{2p\beta }|F|\chi_{\{v_{k} > 0\}})^{\frac{1}{p}}
+ \int_{B(r_{0})}R^{2\beta }|F|log^{+}|F|\cdot \chi_{\{v_{k} > 0 \}} \} ,
\end{split}
\end{equation*}
in which the symbol $(P_{k1})_{B(r_{0})}$ stands for the mean value of $P_{k1}$ over the open ball $B(r_{0})$.
Since we know that $\{v_{k}  > 0\}$ is a subset of $\{|u| > \frac{R}{2}\}$, for all $k\geqslant 1$, so it follows from the above inequality that
\begin{equation*}
\begin{split}
\int_{B(r_{0})}|F|\chi_{\{  v_{k} > 0 \}}|P_{k1} - (P_{k1})_{B(r_{0})}|dx 
&\leqslant \frac{2C_{0}}{\alpha} \frac{p}{p-1}4^{p\beta}
\{ 1 + K^{1-\frac{1}{p}}\}\times\\
&\{ (\int_{B(r_{0})}|u|^{2p\beta}|F|\chi_{\{v_{k} > 0\}})^{\frac{1}{p}}\\
&+ \int_{B(r_{0})} |u|^{2\beta}|F|\log^{+}|F|\cdot \chi_{\{v_{k} > 0\}}\} . 
\end{split}
\end{equation*}
So, we can conclude from inequality (\ref{10}), and the above inequality that
\begin{equation}\label{11}
\begin{split}
\int_{T_{k-1}}^{1}|\int_{\mathbb{R}^{3}}\nabla (\frac{v_{k}}{|u|})u P_{k1}dx|dt
&\leqslant R \frac{2C_{0}}{\alpha }\frac{p}{p-1}4^{p\beta } (1+K^{1-\frac{1}{p}})\times \\ 
&\{(\int_{T_{k-1}}^{1}\int_{B(r_{0})}|u|^{2p\beta }|F|\chi_{\{v_{k} > 0 \}})^{\frac{1}{p}}\\
& + \int_{T_{k-1}}^{1}\int_{B(r_{0})}|u|^{2\beta }|F|\log (1+|F|)\chi_{\{v_{k} > 0 \}} \}\\
& + C(r_{0})2^{\beta }R \|u\|_{L^{\infty}(L^{2})}
\int_{T_{k-1}}^{1}\int_{B(r_{0})}|u|^{\beta }|F|\chi_{\{v_{k} > 0 \}} . 
\end{split}
\end{equation}
In order to use the given hypothesis that $|u\cdot \nabla F| (t,x)   \leqslant A |u(t,x)| |F(t,x)|$, for any $(t,x) \in [0,1)\times B(r_{0})$ satisfying
$|F(t,x)| \geqslant L$ (with $L > 0$ to be the given constant in Theorem \ref{goal}), we carry out the following estimate. 
\begin{equation}\label{12}
\begin{split}
\int_{T_{k-1}}^{1}\int_{B(r_{0})} |u|^{2\beta }|F|\log (1+|F|)\chi_{\{v_{k} > 0\}}
&\leqslant \int_{T_{k-1}}^{1}\int_{B(r_{0})} |u|^{2\beta }|F|\log (1+|F|)\chi_{\{|F|\leqslant L+1\}}\chi_{\{v_{k} > 0   \}}\\
& + \int_{T_{k-1}}^{1}\int_{B(r_{0})}|u|^{2\beta }|F|\log (1+|F|)\chi_{\{|F|> L+1\}}\chi_{\{v_{k} > 0\}}\\
&\leqslant (L+1)\log (L+2) \int_{T_{k-1}}^{1}\int_{B(r_{0})}  |u|^{2\beta }\chi_{\{v_{k} > 0 \}}\\
&+ \int_{T_{k-1}}^{1}\int_{B(r_{0})}  |u|^{2\beta }|F|\log (1+|F|)\chi_{\{|F|> L+1\}}
\chi_{\{v_{k} > 0\}} .
\end{split}
\end{equation}

\noindent{\bf Step five}
To deal with the second term in the last line of inequality (\ref{12}), we consider the sequence $\{\phi_{k} \}_{k=1}^{\infty}$ of nonnegative continuous functions on 
$[0,\infty)$, which are defined by 
\begin{itemize}
\item $\phi_{k} (t)= 0$, for all $t\in [0, C_{k}]$.
\item $\phi_{k} (t)= t - C_{k}$, for all $t\in (C_{k},  C_{k}+1)$.
\item $\phi_{k} (t)= 1$, for all $t\in [C_{k}+1, +\infty)$.
\end{itemize}
where the symbol $C_{k}$ stands for $C_{k} = R(1-\frac{1}{2^{k}})$, for every $k\geqslant 1$. Here, we remark that, for the purpose of taking spatial derivative, the 
composite function $\phi_{k} (|u|)$ is a good substitute for $\chi_{\{v_{k} > 0\}} = \chi_{\{|u| > R (1- \frac{1}{2^{k}})\}}$, since $\phi_{k}$ is Lipschitz. 
Moreover, we also need a smooth function $\psi : \mathbb{R}\rightarrow \mathbb{R}$ satisfying the following conditions that:
\begin{itemize}
\item $\psi (t)=1$, for all $t\geqslant L+1$.
\item $0 < \psi (t) < 1$, for all $t$ with $L< t < L+1$.
\item $\psi (t) = 0$, for all $-L \leqslant t \leqslant L $.
\item $-1 < \psi (t) < 0$, for all $t$ with $-L-1 < t < -L$.
\item $\psi (t) = -1$, for all $t \leqslant -L-1$.
\item $0 \leqslant \frac{d}{dt}\psi  \leqslant 2$, for all $t\in \mathbb{R}$. 
\end{itemize}
We further remark that the smooth function $\psi : \mathbb{R}\rightarrow \mathbb{R}$ characterized by the above properties must also satisfy the property that
$ \psi'(t)=\frac{d \psi }{dt}|_{(t)} = 0  $, on $t \in (-\infty , -L-1) \cup (-L,L) \cup (L+1 , \infty )$, which will be employed in forthcoming inequality estimations \ref{14} 
and \ref{14star} without explicit mention.
With the above preperation, let $\beta$ be such that $\frac{3}{\alpha} < \beta < \frac{10}{3\alpha} $, with $\alpha$ to be the given index as \emph{specified in} Theorem 
\ref{goal}. 

We now consider the function $F = div(\frac{u}{|u|})$, and recall that our solution $u$ satisfies 
$|u\cdot \nabla F|\leqslant A|F|\cdot |u|$ on $\{(t,x) \in [0,1)\times B(r_{0}) : |F(t,x)| \geqslant L \}$. for some given constant $L > 0$.



it follows that
\begin{itemize}
\item $|u\cdot \nabla F|(t,x) \leqslant A (L+1)|u(t,x)|$, if $(t,x) \in [0,1)\times B(r_{0})$ satisfies $ L \leqslant |F(t,x)|\leqslant L+1$.

\item $|\frac{u\cdot \nabla |F|}{1+ |F|}| \leqslant 
\frac{|u\cdot \nabla |F||}{|F|}= \frac{|u\cdot \nabla F|}{|F|} \leqslant A |u|$ is valid on  $[0,1)\times B(r_{0})\cap \{|F(s)|\geq L  \}$.
\end{itemize}
Then, we carry out the following calculation on $[0,1)\times B(r_{0})$, for each $k \geqslant 1 $.
\begin{equation}\label{13}
\begin{split}
div\{ |u|^{2\beta -1} u \psi (F) \log (1+|F|) \phi_{k}(|u|)\}
&= -(2\beta -1)|u|^{2\beta}F \psi (F) \log (1+|F|)\phi_{k}(|u|)\\
&-|u|^{2\beta +1}F\psi (F) \log (1+|F|)\chi_{\{C_{k}  < |u| < C_{k} +1\}}\\
&+ |u|^{2\beta -1} \frac{d\psi }{dt}(F) (u\cdot \nabla F)\log (1+|F|)\phi_{k}(|u|)\\
& + |u|^{2\beta -1}\psi (F) \frac{u\cdot \nabla |F|}{1+|F|} \phi_{k}(|u|)  , 
\end{split}
\end{equation}
Since $R > 2M_{0}$ ensures that, for each $t\in [0,1)$, $\phi_{k}(|u|)(t,\cdot )$ is \emph{compactly supported} in $B(r_{0})$, we have the following equality for each 
$t\in [0,1)$.
\begin{equation*}
\int_{B(r_{0})} div\{ |u|^{2\beta -1} u \psi (F) \log(1+|F|) \phi_{k}(|u|)\} = 0. 
\end{equation*}
So, it follows from inequality (\ref{13}) that
\begin{equation}\label{14}
\begin{split}
\Lambda_{1} + \Lambda_{2}
&\leqslant \int_{T_{k-1}}^{1}\int_{B(r_{0})}|u|^{2\beta -1}|\frac{d\psi }{dt}(F)|\cdot 
|u\cdot \nabla F|\log (1+|F|)\phi_{k}(|u|)\\
&+ \int_{T_{k-1}}^{1}\int_{B(r_{0})}|u|^{2\beta -1}|\psi (F)|\cdot |\frac{u\cdot \nabla |F|}{1+|F|}|
\phi_{k} (|u|)\\
&\leqslant \int_{T_{k-1}}^{1}\int_{B(r_{0})} |u|^{2\beta -1}(2)(A (L+1) |u|)
\log(L+2)\phi_{k}(|u|) \\
&+ \int_{T_{k-1}}^{1}\int_{B(r_{0})} |u|^{2\beta -1}\cdot A\cdot |u|\phi_{k}(|u|) \cdot \chi_{\{|F| \geqslant L \}}\\
&\leqslant A[2(L+1)\log (L+2) + 1] \int_{T_{k-1}}^{1}\int_{B(r_{0})} |u|^{2\beta } \phi_{k}(|u|)\\
&\leqslant A[2(L+1)\log (L+2) + 1] \int_{T_{k-1}}^{1}\int_{B(r_{0})} |u|^{2\beta } \chi_{\{v_{k} > 0 \}} ,
\end{split}
\end{equation}
in which the terms $\Lambda_{1}$, and $\Lambda_{2}$ are given by
\begin{itemize}
\item $\Lambda_{1} = (2\beta -1)\int_{T_{k-1}}^{1}\int_{B(r_{0})} |u|^{2\beta}F\psi (F)\cdot log(1+|F|)
\phi_{k}(|u|)$.
\item $\Lambda_{2} = \int_{T_{k-1}}^{1}\int_{B(r_{0})} |u|^{2\beta +1}(F\psi (F))\cdot log(1+|F|)
\chi_{\{C_{k} < |u| < C_{k} +1 \}}   $ .
\end{itemize}
We then notice that
\begin{itemize}
\item Since $\beta > \frac{3}{\alpha} > 1$, we have $\Lambda_{1} \geqslant
\int_{T_{k-1}}^{1}\int_{B(r_{0})}  |u|^{2\beta }(F\psi (F))\log (1+|F|)\chi_{\{ |u|\geqslant C_{k} +1 \}}   $.
\item $\Lambda_{2} \geqslant \frac{R}{2}  \int_{T_{k-1}}^{1}\int_{B(r_{0})}   |u|^{2\beta}F\psi(F)\log (1+|F|)
\chi_{\{C_{k} < |u| < C_{k} + 1 \}} $, for every $k\geqslant 1$. Notice that this is true because $C_{k}= R(1-\frac{1}{2^k})$, and that 
$(1-\frac{1}{2^{k}}) \geqslant \frac{1}{2}$, for every $k\geqslant 1$.
\end{itemize}
Since $|F|\chi_{\{|F| > L+1\}} \leqslant |F| |\psi(F)| = F\psi(F)$,  it follows from inequality (\ref{14}) that
\begin{equation}\label{15}
\begin{split}
&\int_{T_{k-1}}^{1}\int_{B(r_{0})} |u|^{2\beta }|F| \log (1+|F|)\chi_{\{|F| > L+1\}}\chi_{\{v_{k} > 0\}}\\
& \leqslant \int_{T_{k-1}}^{1}\int_{B(r_{0})} |u|^{2\beta }F\psi (F)\log (1+|F|)\chi_{\{v_{k} > 0\}}\\
&\leqslant  \int_{T_{k-1}}^{1}\int_{B(r_{0})}      |u|^{2\beta }F\psi (F)\log (1+|F|)\chi_{\{C_{k} < |u|
< C_{k} +1 \}}\\
&+ \int_{T_{k-1}}^{1}\int_{B(r_{0})}  |u|^{2\beta }F\psi (F)log(1+|F|)\chi_{\{|u| \geqslant C_{k} +1\}}\\
&\leqslant \frac{2}{R}\Lambda_{2} + \Lambda_{1}\\
&\leqslant 2A[2(L+1)\log(L+2) +1] \int_{Q_{k-1}} |u|^{2\beta}\chi_{\{v_{k} > 0\}}. 
\end{split}
\end{equation}
By using inequality \eqref{weakLPimprove2} in Lemma \ref{weakLP}, we raise up the index for the term $\int_{Q_{k-1}} |u|^{\theta} \chi_{\{v_{k}  > 0\}}$, for any $\theta$ 
with $ 0  < \theta < \frac{10}{3}$, in the following way
\begin{equation*}
\begin{split}
\int_{Q_{k-1}}|u|^{\theta }\chi_{\{v_{k} > 0 \}}
& = \int_{Q_{k-1}} \{R(1-\frac{1}{2^{k}}) + v_{k} \}^{\theta } \chi_{\{v_{k}  > 0 \}}\\
&\leqslant C_{\theta} \{R^{\theta }\int_{Q_{k-1}}\chi_{\{v_{k}  > 0 \}} +
\int_{Q_{k-1}} v_{k}^{\theta }\chi_{\{v_{k} > 0 \}} \}\\
&\leqslant \frac{C_{\theta }}{R^{\frac{10}{3}-\theta }} \{2^{\frac{10k}{3}}
+ 2^{(\frac{10}{3} - \theta  )k}\}\int_{Q_{k-1}}v_{k-1}^{\frac{10}{3}}\\
&\leqslant \frac{C_{\theta }}{R^{\frac{10}{3}-\theta + (\alpha-2)(\frac{2}{3} -\delta ) }}2^{\frac{10k}{3}}  \{ \frac{2^{\alpha -1}}{\alpha - 2} 
\|u\|_{L^{\infty}(L^{\alpha , *})}\}^{\frac{2}{3} - \delta}    U_{k-1}^{1+\delta} ,  
\end{split}
\end{equation*}
for every $\theta$ with $0< \theta < \frac{10}{3}$, where $C_{\theta }$ is some positive constant depending only on $\theta$. Hence it follows from 
inequalities(\ref{12}), (\ref{15}), and our last inequality that
\begin{equation}\label{16}
\begin{split}
\int_{T_{k-1}}^{1}\int_{B(r_{0})} |u|^{2\beta } |F|\cdot log(1+|F|)\chi_{\{v_{k} > 0  \}}
&\leqslant (L+1)\log (L+2) \int_{T_{k-1}}^{1}\int_{B(r_{0})} |u|^{2\beta }\chi_{\{v_{k} > 0\}}\\
&+ \int_{T_{k-1}}^{1}\int_{B(r_{0})} |u|^{2\beta }|F|log(1+|F|)\chi_{\{|F|> L+1\}}
\chi_{\{v_{k} > 0 \}}\\
&\leqslant \frac{(L+1)\log (L+2) C_{2\beta }2^{\frac{10k}{3}}}{R^{\frac{10}{3}-2\beta +(\alpha-2)(\frac{2}{3} - \delta ) }}\{ \frac{2^{\alpha -1}}{\alpha - 2} 
\|u\|_{L^{\infty}(L^{\alpha , *})}\}^{\frac{2}{3} - \delta}    U_{k-1}^{1+\delta}\\
&+ C_{(A,L)} \int_{Q_{k-1}}|u|^{2\beta }\chi_{\{v_{k} > 0\}}\\
&\leqslant C_{(\beta , A,L)}\cdot 2^{\frac{10k}{3}} \{ \frac{2^{\alpha -1}}{\alpha - 2} \|u\|_{L^{\infty}(L^{\alpha , *})}\}^{\frac{2}{3} - \delta}  U_{k-1}^{1+\delta } \\
& \times \{ \frac{1}{R^{\frac{10}{3}-2\beta +  (\alpha-2)(\frac{2}{3}-\delta )   }}  \} ,
\end{split}
\end{equation}

in which $\beta > \frac{3}{\alpha }$, and that $\beta$ is sufficiently close to $\frac{3}{\alpha}$, and $C_{\beta , A,L }$ is some constant depending only on $\beta$, 
$A$, and $L$. 
Next, we also need to deal with  
$( \int_{T_{k-1}}^{1}\int_{B(r_{0})} |u|^{2p\beta }|F|\chi_{\{v_{k}\geqslant 0 \}} )^{\frac{1}{p}}$,
and $ \int_{T_{k-1}}^{1}\int_{B(r_{0})} |u|^{\beta }|F|\chi_{\{v_{k}\geqslant 0\}}$,
which appear in inequality (\ref{11}). For this purpose, we will consider $\lambda$ which satisfies $\frac{3}{\alpha} < \lambda < \frac{10}{3}$ (we will take $\lambda$ to be 
$2p\beta$ and $\beta$ respectively in forthcoming inequality estimates \ref{17} and \ref{18} ), and let us carry out the following computation, in which $\psi$ and 
$\phi_{k}$ etc are just the same as before. 
\begin{equation*}
\begin{split}
div \{|u|^{\lambda -1 }u\psi (F) \phi_{k}(|u|)\}
&= -(\lambda -1)|u|^{\lambda }F\psi (F)\phi_{k}(|u|)\\
& + |u|^{\lambda -1}\frac{d\psi }{dt}(F) (u\cdot \nabla F )\phi_{k}(|u|)\\
&- |u|^{\lambda +1}F\psi (F)\chi_{\{C_{k} < |u| < C_{k} +1   \}} .
\end{split}
\end{equation*}
Since $R > 2M_{0}$ ensures that $\phi_{k}(|u|)$ is compactly supported in $B(r_{0})$, we have, for each $t\in [0,1)$, that
\begin{equation*}
\int_{B(r_{0})} div \{|u|^{\lambda -1 }u\psi (F) \phi_{k}(|u|)\} = 0 .
\end{equation*}
Hence, it follows from $|\frac{d\psi }{dt}(F)| \leqslant 2 \chi_{\{L < |F| < L+1 \}}$ and the above equality that
\begin{equation}\label{14star}
\begin{split}
(\lambda -1 ) \int_{T_{k-1}}^{1}\int_{B(r_{0})} |u|^{\lambda}F\psi (F)\phi_{k} (|u|)
& + \int_{T_{k-1}}^{1}\int_{B(r_{0})} |u|^{\lambda +1}F\psi (F)\chi_{\{C_{k} < |u| < C_{k} +1 \}}\\
&\leqslant \int_{T_{k-1}}^{1}\int_{B(r_{0})} |u|^{\lambda -1}|\frac{d\psi }{dt}(F)|\cdot 
|u\cdot \nabla F|\phi_{k}(|u|)\\
&\leqslant \int_{Q_{k-1}} |u|^{\lambda -1} (2)(A (L+1)|u|)
\chi_{\{v_{k}  > 0\}}\\
&\leqslant 2A (L+1)\int_{Q_{k-1}} |u|^{\lambda }\chi_{\{v_{k} > 0 \}} .
\end{split}
\end{equation}
By the same calculation as in inequality (\ref{14}), we can see that
\begin{equation*}
\begin{split}
 \int_{T_{k-1}}^{1}\int_{B(r_{0})} |u|^{\lambda } F\psi (F) \chi_{\{v_{k}  > 0 \}}
& \leqslant  \int_{T_{k-1}}^{1}\int_{B(r_{0})} |u|^{\lambda }F\psi (F)\chi_{\{C_{k} < |u| < C_{k} + 1   \}}\\
&+  \int_{T_{k-1}}^{1}\int_{B(r_{0})} |u|^{\lambda } F\psi (F)\chi_{\{|u| \geqslant C_{k} + 1  \}}\\
&\leqslant \frac{2}{R}  \int_{T_{k-1}}^{1}\int_{B(r_{0})}  |u|^{\lambda + 1} F\psi (F)
\chi_{\{C_{k} < |u| < C_{k} + 1 \}}\\
& +  \int_{T_{k-1}}^{1}\int_{B(r_{0})} |u|^{\lambda }F\psi (F) \phi_{k}(|u|)\\
&\leqslant (2 + \frac{1}{\lambda - 1} ) \{   \int_{T_{k-1}}^{1}\int_{B(r_{0})}  |u|^{\lambda +1}F\psi (F)\chi_{\{C_{k}  < |u| < C_{k} +1 \}}\\    
& +(\lambda -1)  \int_{T_{k-1}}^{1}\int_{B(r_{0})}   |u|^{\lambda }F\psi(F)\phi_{k}(|u|)\}\\
&\leqslant 2A (L+1) (2+\frac{1}{\lambda - 1}) \int_{Q_{k-1}}|u|^{\lambda }\chi_{\{v_{k} > 0 \}},
\end{split}
\end{equation*}
in which $\lambda$ satisfies $\frac{3}{\alpha} < \lambda < \frac{10}{3} $. Now, put $\lambda = 2p\beta $, with $\beta > \frac{3}{\alpha}$ to be sufficiently close to 
$\frac{3}{\alpha}$, and $p>1$ to be sufficiently close to $1$. Since $|F|\chi_{\{|F|> L+1\}} \leqslant |F||\psi(F)| = F \psi(F)$,
 it follows from our last inequality that
\begin{equation}\label{17}
\begin{split}
 \int_{T_{k-1}}^{1}\int_{B(r_{0})} |u|^{2p\beta }|F|\chi_{\{v_{k} > 0 \}}
& =  \int_{T_{k-1}}^{1}\int_{B(r_{0})} |u|^{2p\beta }|F|\chi_{\{|F|\leqslant L+1\}}
\chi_{\{v_{k} > 0\}}\\
&  +   \int_{T_{k-1}}^{1}\int_{B(r_{0})}  |u|^{2p\beta }\chi_{\{|F|> L+1\}}
\chi_{\{v_{k} > 0  \}}|F|\\
&\leqslant (L+1)\int_{Q_{k-1}}|u|^{2p\beta }\chi_{\{v_{k}  > 0\}}\\
&+ 2A (L+1) (2+\frac{1}{2p\beta - 1})\int_{Q_{k-1}}|u|^{2p\beta }\chi_{\{v_{k} > 0 \}}\\
&\leqslant  \frac{C_{(\beta , A,L,p)}}{R^{\frac{10}{3}-2p\beta + (\alpha-2)(\frac{2}{3} -\delta ) }}   
\cdot  2^{\frac{10k}{3}} \{\frac{2^{\alpha-1}}{\alpha -2} \|u\|_{L^{\infty}(L^{\alpha ,*})}\}^{\frac{2}{3} - \delta } U_{k-1}^{1+\delta} .
\end{split}
\end{equation}
In exactly the same way, by setting $\lambda$ to be $\beta $, with 
$\beta > \frac{3}{\alpha}$ to be sufficiently close to $\frac{3}{\alpha}$, it also follows that
\begin{equation}\label{18}
\begin{split}
\int_{T_{k-1}}^{1}\int_{B(r_{0})}  |F|\chi_{\{v_{k} > 0 \}}
&= \int_{T_{k-1}}^{1}\int_{B(r_{0})}  |u|^{\beta }|F|\chi_{\{|F|\leqslant L+1\}}\chi_{\{v_{k} > 0\}}\\
&+   \int_{T_{k-1}}^{1}\int_{B(r_{0})}  |u|^{\beta }|F|\chi_{\{|F|> L+1\}}\chi_{\{v_{k} > 0 \}}\\
&\leqslant (L+1)\int_{Q_{k-1}}|u|^{\beta }\chi_{\{v_{k} > 0 \}}
+2A(L+1)(2+\frac{1}{\beta - 1})\int_{Q_{k-1}}|u|^{\beta }\chi_{\{v_{k} > 0 \}}\\
& \leqslant  \frac{C_{(\beta ,A, L)}}{R^{\frac{10}{3}- \beta + (\alpha-2)(\frac{2}{3} -\delta ) }}   
\cdot  2^{\frac{10k}{3}} \{\frac{2^{\alpha-1}}{\alpha -2} \|u\|_{L^{\infty}(L^{\alpha ,*})}\}^{\frac{2}{3} - \delta } U_{k-1}^{1+\delta} .
\end{split}
\end{equation}
By combining inequalities (\ref{11}), (\ref{16}), and (\ref{17}),and (\ref{18}) we now conclude that
\begin{equation}\label{19}
\begin{split}
\int_{Q_{k-1}}|\int_{Q_{k-1}} \nabla (\frac{v_{k}}{|u|})uP_{k1}dx|ds
&\leqslant (1+\frac{1}{\alpha })C(\beta ,A, L, p )(1+K^{1-\frac{1}{p}}) (1+\|u\|_{L^{\infty}(L^{2}))}) \\
& \{ [\frac{2^{\alpha-1}}{\alpha -2} \|u\|_{L^{\infty}(L^{\alpha ,*})}]^{\frac{2}{3} - \delta } + 
 [\frac{2^{\alpha-1}}{\alpha -2} \|u\|_{L^{\infty}(L^{\alpha ,*})}]^{(\frac{2}{3} - \delta)\frac{1}{p} } \}\\
&\{  (\frac{1}{R^{\frac{10}{3}-2p\beta +(\alpha-2)(\frac{2}{3} - \delta ) - p}})^{\frac{1}{p}} 2^{\frac{10k}{3p}} U_{k-1}^{\frac{1}{p}(1+\delta )} \\  
& + \frac{1}{R^{\frac{10}{3}-2\beta + (\alpha-2)(\frac{2}{3} -\delta ) -1 }} 2^{\frac{10k}{3}}U_{k-1}^{1+\delta }
 \} .
\end{split}
\end{equation}
Before we proceed to the last step and complete the proof of Theorem \ref{goal}, let us briefly explain why the condition $1 + 2(\frac{\alpha}{3} -\frac{3}{\alpha}) > 0$
imposed on $2< \alpha < 3$ is necessary. Notice that if $p\rightarrow 1^{+}$, and $\beta \rightarrow \frac{3}{\alpha}^{+}$, and $\delta \rightarrow 0^{+}$ , we have 
$(\frac{10}{3}-2p\beta +(\alpha-2)(\frac{2}{3} - \delta ) - p )\rightarrow 1 + 2(\frac{\alpha}{3} -\frac{3}{\alpha}) $, and that 
$(\frac{10}{3}-2\beta + (\alpha-2)(\frac{2}{3} -\delta ) -1 )\rightarrow  1 + 2(\frac{\alpha}{3} -\frac{3}{\alpha}) $.
This explains that the condition $1 + 2(\frac{\alpha}{3} -\frac{3}{\alpha}) > 0$ on $\alpha \in (2,3)$ is necessary if we insist that
both $(\frac{10}{3}-2p\beta +(\alpha-2)(\frac{2}{3} - \delta ) - p )$ and $(\frac{10}{3}-2\beta + (\alpha-2)(\frac{2}{3} -\delta ) -1 ) $ have to be positive.\\

\noindent{\bf Step Six: Final step of the proof}

By combining inequalities \eqref{trival}, \eqref{Divide}, (\ref{middle}), (\ref{easy}), and (\ref{19}), we conclude that the following estimate is valid.

\begin{equation}\label{Final}
\begin{split}
U_{k} & \leqslant \frac{2^{\frac{10k}{3}}}{R^{\frac{4}{3}}}C_{0}U_{k-1}^{\frac{5}{3}} +
 C(\beta ,A ,L , p, \delta , \|u\|_{L^{\infty}L^{2}} ,\|u\|_{L^{\infty}L^{\alpha ,*}}  ) \{\frac{  U_{k-1}^{\frac{1}{p} + \delta (\frac{2-p}{2p})}  }{R^{\beta [\frac{10-8p}{3p} + (\frac{2-p}{2p})(\alpha - 2) (\frac{2}{3} -\delta )] -(\frac{2-p}{p})}} \\
 & +  (\frac{  U_{k-1}^{\frac{1}{p}(1+\delta )} }{R^{\frac{10}{3}-2p\beta +(\alpha-2)(\frac{2}{3} - \delta ) - p}})^{\frac{1}{p}}  
+ \frac{U_{k-1}^{1+\delta }}{R^{\frac{10}{3}-2\beta + (\alpha-2)(\frac{2}{3} -\delta ) -1 }} \}
\end{split}
\end{equation}

Here, in order to derive the conclusion $|u| \leqslant [\frac{3}{4} , 1)\times \mathbb{R}^{3}$ by using inequality \eqref{Final}, we have to be very careful in the 
selection of the constants $\beta$, $p$, $\delta$ etc. This is due to the following fact. On the one hand, we require all the powers of $U_{k-1}$ such as 
$\frac{1}{p} + \delta (\frac{2-p}{2p})$, $\frac{1}{p}(1+\delta)$, and $1+\delta$ to be strictly positive, so that $p$ has to be sufficiently close to $1$ and 
that $\delta$, however small, has to stay positive. On the other hand, the constant 
$C(\beta ,A ,L , p, \delta , \|u\|_{L^{\infty}L^{2}} , \|u\|_{L^{\infty}L^{\alpha ,*}}  )$ will blow up to $\infty$
if $p \rightarrow 1^{+}$. So, to clarify the situation, we have to fix the choice of $\beta$ \emph{first} by using the condition
$1 + 2(\frac{\alpha}{3} -\frac{3}{\alpha}) > 0$ on $\alpha \in (2,3)$. Once the choice of $\beta$ is fixed, we will \emph{fix} the parameters $p > 1$ and $\delta > 0$.\

Observe that the condition $1 + 2(\frac{\alpha}{3} -\frac{3}{\alpha}) > 0$ on  $\alpha \in (2,3)$ is equivalent to $\frac{1}{2} + \frac{\alpha}{3} > \frac{3}{\alpha}$, and
this allows us to select some $\beta$ to be in the interval $(\frac{3}{\alpha} , \frac{1}{2} + \frac{\alpha}{3}  )$. Now, let $\beta$ to be a \emph{fixed} choice of
positive number which satisfies $\frac{3}{\alpha} < \beta < \frac{1}{2} + \frac{\alpha}{3} $. Next, recall that we have the following limiting relations.

\begin{itemize}
\item $ \lim_{p\rightarrow 1^{+} ,\delta \rightarrow 0^{+} } \beta [\frac{10-8p}{3p} + (\frac{2-p}{2p})(\alpha - 2) (\frac{2}{3} -\delta )] -(\frac{2-p}{p}) = 
\beta(\frac{\alpha}{3}) - 1.    $
\item $\lim_{p\rightarrow 1^{+} ,\delta \rightarrow 0^{+} }  \{\frac{10}{3}-2p\beta +(\alpha-2)(\frac{2}{3} - \delta ) - p \} = 
2\{\frac{1}{2} + \frac{\alpha}{3}  -\beta    \}           $.
\item $ \lim_{\delta \rightarrow 0^{+}} \frac{10}{3}-2\beta + (\alpha-2)(\frac{2}{3} -\delta ) -1 = 2\{\frac{1}{2} + \frac{\alpha}{3}  -\beta    \}    $ .
\end{itemize}

Notice that the fixed choice of $\beta$ with $\frac{3}{\alpha} < \beta < \frac{1}{2} + \frac{\alpha}{3} $ ensures that the limiting constants $\beta(\frac{\alpha}{3}) - 1 $
and $ 2\{\frac{1}{2} + \frac{\alpha}{3}  -\beta    \}   $ are both positive simultaneously. As a result, the above three limiting relations imply that for some 
\emph{fixed} choice of $p > 1$ sufficiently close to $1$, and some \emph{fixed} choice of $\delta > 0$ sufficiently close to $0$ 
(both depending on the choice of $\beta$), it follows that
the following three constants are \emph{positive}.

\begin{itemize}
\item  $\beta [\frac{10-8p}{3p} + (\frac{2-p}{2p})(\alpha - 2) (\frac{2}{3} -\delta )] -(\frac{2-p}{p}) > 0 $.
\item  $\{\frac{10}{3}-2p\beta +(\alpha-2)(\frac{2}{3} - \delta ) - p \} > 0$.
\item  $\frac{10}{3}-2\beta + (\alpha-2)(\frac{2}{3} -\delta ) -1 > 0$.
\end{itemize}

This observation allows us to use nonlinear recurrence relation \eqref{Final} to deduce that as long as $R > M_{0} + 1$
is chosen to be sufficiently large, $U_{1}$ will become smaller than the universal constant $C_{0}^{*}$ as required by Lemma \ref{Vass}. According to Lemma \ref{Vass}, 
this smallness of $U_{1}$ will lead to the decay of $U_{k}$ to $0$ as $k \rightarrow \infty$, and this in turn will lead to the conclusion that $|u| \leqslant R$ is valid
over $[\frac{3}{4} , 1) \times \mathbb{R}^{3}$, for some sufficiently large constant $R$. Hence, it follows that the smoothness of $u$ can be extended beyond the 
possible blow up time $1$.\\

\textbf{Acknowledgments: }Both authors are grateful to Professor Vladimir {\v{S}}ver{\'a}k for his encouragement and guidance.
This paper was developed during
a stay of the second author at the Institute for Mathematics and Its Applications, University of Minnesota. 

\bibliography{largeswirl.bib}

\end{document}